\documentclass[a4paper,12pt, reqno]{amsart}
\usepackage{amssymb,amsthm,amsmath}
\usepackage{cite}

\usepackage{amsmath,amsthm,amscd,amssymb}
\usepackage{latexsym}
\usepackage[colorlinks,citecolor=red,pagebackref,hypertexnames=false]{hyperref}

\numberwithin{equation}{section}
\allowdisplaybreaks[2]
\theoremstyle{plain}
\newtheorem{theorem}{Theorem}[section]

\newtheorem{lemma}[theorem]{Lemma}
\newtheorem{corollary}[theorem]{Corollary}
\newtheorem{proposition}[theorem]{Proposition}

\theoremstyle{definition}
\newtheorem{definition}[theorem]{Definition}

\newtheorem{example}[theorem]{Example}

\theoremstyle{remark}
\newtheorem{remark}[theorem]{Remark}

\newtheorem{case[theorem]}{Case}

\def\norm#1.#2.{\lVert#1\rVert_{#2}}

\title[Szeg\"o type limit theorems on the Heisenberg group]
{Szeg\"o type limit theorems on the Heisenberg group}

\author{Shyam Swarup Mondal}
\author{Jitendriya Swain}
\address{Department of Mathematics,
	Indian Institute of Technology Guwahati,
	Guwahati 781039, \;\; India.} 
\email{shyam.mondal@iitg.ac.in, jitumath@iitg.ac.in}

\begin{document}
	\allowdisplaybreaks
\baselineskip=18.5pt

\begin{abstract} Let $\mathcal{H}=-\Delta_{\mathbb{H}}+V$ be the Schr\"odinger operator on the Heisenberg group $\mathbb{H}^n$, where $\Delta_{\mathbb{H}}$ is the full laplacian on $\mathbb{H}^n$ and  $V$ is  a positive smooth potential, bounded below and grows like $|g|^\kappa, \kappa>0$ for large $|g|$.  Let $\mathcal{P}_{r}$ be the orthogonal projection of $L^2(\mathbb{H}^n)$ onto the space of eigenfunctions of $\mathcal{H}$ with eigenvalue $\leq r$; Let $\textbf{b}$ be a bounded real valued integrable function on $\mathbb{H}^n$ and $M_{\textbf{b}}$ be the operator of multiplication by $\textbf{b}$ on $L^2(\mathbb{H}^n)$.  Then for any $f \in C(\mathbb{R})$ we have
\begin{align*}
\lim_{r\to\infty} \frac{{tr}{f(\mathcal{P}_rM_{\textbf{b}}\mathcal{P}_r)}}{tr~(\mathcal{P}_r)} &=  \int_{\mathbb{H}^n}f(\textbf{b}(g)) \,dg.\end{align*} Further, if $A$ be a 0-th order self-adjoint pseudo-differential operator on $L^2(\mathbb{H}^n)$ relative to the operator      $1+|\lambda|H+V(g), g\in \mathbb{H}^n, \lambda \in \mathbb{R}^*$ with symbol $a(g, {\lambda}),$ where $H$ is the Hermite operator on $L^2(\mathbb{R}^n)$  then
	\begin{align*}
	\lim_{r\to\infty} \frac{tr~{f(\mathcal{P}_rA\mathcal{P}_r)}}{tr~(\mathcal{P}_r)} &= \lim_{r\to\infty}  \frac{\int_{G^{r}}f(a_{g, {\lambda}}(\xi, x))  \,d\xi\,dx \,dg\,d\mu(\lambda) }{\int_{G^{r}}  \,d\xi\,dx \,dg\,d\mu(\lambda)},
	\end{align*}
	(Assuming one limit exists)\\
	where $G^{r}=\{(g, \lambda, \xi, x)\in \mathbb{H}^n \times \mathbb{R}^*\times  \mathbb{R}^n\times \mathbb{R}^n : |\lambda |(1+|\xi| ^2+|x|^2)+V(g)\leq r \}$, $a(g, {\lambda})=Op^W(a_{g, {\lambda}})$, and $\mu(\lambda)$ is the Plancherel measure on the Heisenberg group. Also we show that the above limit on the right hand side remains unaltered under a compact perturbation of the pseudo-differential operator $A$ or a perturbation of the Schr\"odinger operator $\mathcal{H}$ by bounded self-adjoint operators on $L^2(\mathbb{H}^n)$.

%
%

\end{abstract}

\date{\today}
\maketitle
\def\BC{{\mathbb C}} \def\BQ{{\mathbb Q}}
\def\BR{{\mathbb R}} \def\BI{{\mathbb I}}
\def\BZ{{\mathbb Z}} \def\BD{{\mathbb D}}
\def\BP{{\mathbb P}} \def\BB{{\mathbb B}}
\def\BS{{\mathbb S}} \def\BH{{\mathbb H}}
\def\BE{{\mathbb E}}
\def\BN{{\mathbb N}}
\def\LP{{W(L^p(\BR^d, \BH), L^q_v)}}
\def\LPN{{W_{\BH}(L^p, L^q_v)}}
\def\LPQ{{W_{\BH}(L^{p'}, L^{q'}_{1/v})}}
\def\L1{{W_{\BH}(L^{\infty}, L^1_w)}}
\def\LB{{L^p(Q_{1/ \beta}, \BH)}}
\def\SP{S^{p,q}_{\tilde{v}}(\BH)}
\def\f{{\bf f}}
\def\h{{\bf h}}
\def\hp{{\bf h'}}
\def\m{{\bf m}}
\def\g{{\bf g}}
\def\ga{{\boldsymbol{\gamma}}}

\section{Introduction}
The observable quantities in the classical system are described by real valued functions on the phase space whereas in quantum systems they are given by self-adjoint operators on a Hilbert space. Therefore it is important to study the correspondence between the classical and quantum statistical mechanics.
Pseudo-differential operator theory provides a natural platform to relate the classical and quantum mechanics. For instance in \cite{zel83}, Zelditch considered the Schr\"{o}dinger operator on $\mathbb{R}^n$ of the form $\widetilde{H}=-\frac{1}{2}\Delta+V$, where $V$ is a smooth
positive function that grows like $V_0|x|^\kappa ,~\kappa>0$ at infinity.
He took a  $0$-th
order self-adjoint pseudo-differential operator $A$ associated with a
symbol $a(x,\xi)$
relative to Beals-Fefferman weights
$\varphi_1(x,\xi)=1,\varphi_2(x,\xi) =(1+|\xi|^2+V(x))^{1/2}$ and proved the
following Szeg\"o type theorem:
For any continuous function $f$,
\begin{align}\label{Newnumber}
\lim_{\lambda \rightarrow \infty} \frac{{ tr} f(P_\lambda A P_\lambda)}{\mbox{rank}~ (P_\lambda)} =
\lim_{\lambda \rightarrow \infty} \frac{\int_{\widetilde{H}(x,\xi) \leq \lambda} f(a(x,\xi)) ~ dxd\xi}{\displaystyle\mbox{Vol}{(\widetilde{H}(x,\xi) \leq \lambda)} },	
\end{align}
where $\widetilde{H}(x,\xi)=\frac{1}{2}|\xi|^2+V(x)$ and $P_{\lambda}$ is the orthogonal projection of $L^2(\mathbb{R}^n)$ onto the space of the eigenfunctions of $\widetilde{H}$ with eigenvalue less equal to $\lambda $,  assuming one limit exists.
Such asymptotic spectral formulae expressing the relation between functions of
pseudo-differential operators and their symbols is an important and
interesting problem in mathematical analysis. We refer to \cite{gui, hor, JZ, sm, wim} for similar results in the literature.

We consider  the Schr\"odinger operator $\mathcal{H}=-\Delta_{\mathbb{H}}+V$  on the Heisenberg group $\mathbb{H}^n$, where $\Delta_{\mathbb{H}}$ is the full laplacian on $\mathbb{H}^n$ and  $V$ is  a positive smooth potential, bounded below and grows like $|g|^\kappa, \kappa>0$ for large
  \begin{eqnarray}\label{g}|g|:=(|x|^4+|t|^2)^{\frac{1}{4}},\quad g=(x,t)\in \mathbb{H}^n,\end{eqnarray} defining the homogenous norm on $\mathbb{H}^n$. Such operators are well known to have purely discrete spectrum whose eigenfunctions form a complete set orthonormal basis for $L^2(\mathbb{H}^n)$ (see Theorem 2 of \cite{si} and the $L^2-L^\infty$ boundedness of $e^{-t\Delta_{\mathbb{H}}}$ can be obtained  from (2.2.1) of \cite{kr}). Let $A=Op(Op^W(a_{g, \lambda}))$  be a bounded self-adjoint 0-th order pseudo-differential operator on $L^2(\mathbb{H}^n)$ relative  to the operator   $1+|\lambda|H+V(g)$ (defined in Subsection \ref{subsection1}). For each $r>0$, $\mathcal{P}_rA\mathcal{P}_r$ is a finite rank symmetric operator with spectral measure defined as the sum of  Dirac delta  functions at its eigen values. We show that the sequence of measures $ \frac{{tr}{f(\mathcal{P}_r A\mathcal{P}_r)}}{tr~(\mathcal{P}_r)}$ converges to the weak limit $ \frac{\int_{G^{r}}f(a_{g, {\lambda}}(\xi, x))  \,d\xi\,dx \,dg\,d\mu(\lambda) }{\int_{G^{r}}  \,d\xi\,dx \,dg\,d\mu(\lambda)}$. In particular, if $\textbf{b}$ is a bounded real valued integrable function on $\mathbb{H}^n$ then we obtain the following result with respect to the operator of multiplication 
  $M_{\textbf{b}}$:

\begin{theorem}\label{mo}
Consider the Schr\"odinger operator of the form  $\mathcal{H}=-\Delta_{\mathbb{H}}+V$ on the Heisenberg group $\mathbb{H}^n$.  Let $\mathcal{P}_r$ be the orthogonal projection of $L^2(\mathbb{H}^n)$ onto the space of eigenfunctions of $\mathcal{H}$ with eigenvalue  $\leq r.$ Let $\textbf{b}$ be a bounded real valued integrable function on $\mathbb{H}^n$ and $M_{\textbf{b}}$ be the operator of multiplication by $\textbf{b}$ on $L^2(\mathbb{H}^n)$.  Then for any $f \in C(\mathbb{R})$ we have
\begin{align*}
\lim_{r\to\infty} \frac{{tr}{f(\mathcal{P}_r M_{\textbf{b}}\mathcal{P}_r)}}{tr~(\mathcal{P}_r)} &=  \int_{\mathbb{H}^n}f(\textbf{b}(g)) \,dg.
\end{align*}
\end{theorem}

 We generalize Theorem \ref{mo} by taking a 0-th order self-adjoint pseudo-differential operator on $L^2(\mathbb{H}^n)$ relative to the operator $1+|\lambda|H+V(g)$, where $H$ is the Hermite operator on $L^2(\mathbb{R}^n)$ and  $\lambda\in \mathbb{R}^* $, in place of the multiplication operator $M_{\textbf{b}}$ and  obtain the following Szeg\"o type limit theorem:
\begin{theorem}\label{sch}
Consider the Schr\"odinger operator of the form  $\mathcal{H}=-\Delta_{\mathbb{H}}+V$ on the Heisenberg group $\mathbb{H}^n$. Let $\mathcal{P}_r$ be the orthogonal projection of $L^2(\mathbb{H}^n)$ onto the space of eigenfunctions of $\mathcal{H}$ with eigenvalue $\leq r$; let $A$ be a 0-th order self-adjoint pseudo-differential operator relative to the operator $1+|\lambda|H+V(g)$ on $L^2(\mathbb{H}^n)$ with symbol $a(g, {\lambda})$, where $g\in \mathbb{H}^n, \lambda \in \mathbb{R}^*$ and let $f \in C(\mathbb{R})$. Then
	\begin{align}\label{lim}
	\lim_{r\to\infty} \frac{tr~{f(\mathcal{P}_rA\mathcal{P}_r)}}{tr~(\mathcal{P}_r)} &= \lim_{r\to\infty}  \frac{\int_{G^{r}}f(a_{g, {\lambda}}(\xi, x))  \,d\xi\,dx \,dg\,d\mu(\lambda) }{\int_{G^{r}}  \,d\xi\,dx \,dg\,d\mu(\lambda)},
	\end{align}
	(Assuming one limit exists)\\
	where $G^{r}=\{(g, \lambda, \xi, x)\in \mathbb{H}^n \times \mathbb{R}^*\times  \mathbb{R}^n\times \mathbb{R}^n : |\lambda |(1+|\xi| ^2+|x|^2)+V(g)\leq  r\}$, $a(g, {\lambda})=Op^W(a_{g, {\lambda}})$, and $\mu(\lambda)$ is the Plancherel measure on the Heisenberg group.
\end{theorem}

We also show that the right hand limit in (\ref{lim}) remains unaltered under a perturbation of the Schr\"odinger operator by a bounded self-adjoint operator  $B$ on $L^2(\mathbb{H}^n)$ such that $B+\mathcal{H}$ has discrete spectrum and  
the eigenfunctions of $B+\mathcal{H} $ form a complete orthogonal basis for $L^2(\mathbb{H}^n)$. Note that the operator  $e^{-t(B+\mathcal{H})}=e^{-tB}e^{-t\mathcal{H}}$ is a  compact operator   as 
$e^{-tB}$ is a bounded operator for any $t>0$ (see Theorem 2 of \cite{si}).

\begin{theorem}\label{bdd}
Consider the operator $\mathcal{H}_1=B+\mathcal{H}$   on the Heisenberg group $\mathbb{H}^n$,  where $B$ is a bounded self-adjoint operator on $\mathbb{H}^n$ such that $\mathcal{H}_1$  has purely discrete spectrum and the eigenfunctions of $\mathcal{H}_1$ form a complete orthogonal basis for $L^2(\mathbb{H}^n)$.
Let $\mathcal{P}_r^{'}$ be the orthogonal projection of $L^2(\mathbb{H}^n)$ onto the space of eigenfunctions of $\mathcal{H}_1$ with eigenvalue $\leq r$; let $A$ be a 0-th order self-adjoint pseudo-differential operator relative to the operator $1+|\lambda|H+V(g)$ on $L^2(\mathbb{H}^n)$  with symbol $a(g, {\lambda})$, where $g\in \mathbb{H}^n, \lambda \in \mathbb{R}^*$ and let $f \in C(\mathbb{R})$. Then
\begin{align*}
\lim_{r\to\infty} \frac{{tr}{f(\mathcal{P}_r^{'} A\mathcal{P}_r^{'})}}{tr~(\mathcal{P}_r^{'})}  &= \lim_{r\to\infty}  \frac{\int_{G^{r}}f(a_{g, {\lambda}}(\xi, x))  \,d\xi\,dx \,dg\,d\mu(\lambda) }{\int_{G^{r}}  \,d\xi\,dx \,dg\,d\mu(\lambda)}
	\end{align*}
	(Assuming one limit exists)\\
	where $G^{r}=\{(g, \lambda, \xi, x)\in \mathbb{H}^n \times \mathbb{R}^*\times  \mathbb{R}^n\times \mathbb{R}^n : |\lambda |(1+|\xi| ^2+|x|^2)+V(g)\leq r \}$, $a(g, {\lambda})=Op^W(a_{g, {\lambda}})$, and $\mu(\lambda)$ is the Plancherel measure on the Heisenberg group.
\end{theorem}
We show that the above theorems are valid under a compact perturbation of the pseudo-differential operator $A$ in  Corollary \ref{cpt}.   

To establish our main results we need to consider the ratios of distributions associated to different measures and their asymptotic behaviours. The asymptotic limit of such ratios is compted using Tauberian theorem.  For instance,  Zelditch  in  \cite{zel83}   used  Karamata's Tauberian theorem (see \cite{wid41}), whereas Robert  \cite{rob83} used   the Keldysh Tauberian theorem (see \cite{kel51}).  
However, we use the recent version of Tauberian theorem of  Keldysh by   Grishin-Poedintseva  \cite{gri} 
and a theorem of  Laptev-Safarov   \cite{LapSaff, lap91} for estimate the error term to prove our main results.

Also  we  provide an alternative proof of  the error estimate for $\kappa\in (0,1)$ without using pseudo-differential  symbolic calculus, but  by  proving the boundedness of the operators $[A, V]$ and $[A, \mathcal{L}]$ on $L^2(\mathbb{H}^n)$.

 We   build up the calculus of symbols the pseudo-differential operators relative to the operator $1+|\lambda|H+V(g)$ on $L^2(\mathbb{H}^n)$ using similar techniques used in \cite{ruz2014,ruz14} and establish the link between these symbols and the scalar valued $(\lambda,V(g))$-Shubin classes. Then we construct pseudo-differential approximations to the operator $(\mathcal{H}+u)^{-m}$ on $L^2(\mathbb{H}^n)$ and  $(1+|\lambda|(H+I)+V(g)+u)^{-m}$ on $L^2(\mathbb{R}^n)$ within the calculus of symbols defined related to  $1+|\lambda|H+V(g)$ and $1+ |\lambda |(1+|\xi| ^2+|x|^2)+V(g)$ respectively. Constructing pseudo-differential approximations is almost classical. We refer to \cite{atiyah}  for a detailed study.

We organize the paper as follows.
In Section \ref{sec2},  we provide necessary background  on the Hermite operator, pseudo-differential operators on $\mathbb{R}^n$, and discuss some basic results on the Heisenberg group. In Section \ref{sec3}, we develop the calculus of symbols relative to the operator $1+|\lambda|H+V(g)$ on $L^2(\mathbb{H}^n)$ and establish the link between these symbols and the scaler valued $(\lambda, V(g))$-Shubin class symbols. We construct pseudo-differential approximation to the operator $(\mathcal{ H}+u)^{-m} $ on $L^2(\mathbb{H}^n)$  in Section \ref{sec4}. In Section \ref{sec5} and \ref{sec6}, we prove our main results Theorem \ref{mo}, \ref{sch}, and  \ref{bdd}. Finally, we show that our main results are valid under a compact perturbation of the pseudo-differential operator $A$. We conclude with an alternative proof of the error estimate without using pseudo-differential calculus for $\kappa \in (0, 1)$
\section{Notations and Background}\label{sec2}
The main aim of this section is to define  the symbol classes on the Heisenberg group via the left invariant vector fields and their correspondence with certain symbol classes on $\mathbb{R}^n.$ We start with the definition of the Hermite operator.
\subsection{Hermite Operator}
Let $H_k$ denote the Hermite polynomial on $\mathbb{R}$, defined by
$$H_k(x)=(-1)^k \frac{d^k}{dx^k}(e^{-x^2} )e^{x^2}, \quad k=0, 1, 2, \cdots   ,$$\vspace{.30mm}
and $h_k$ denote the normalized Hermite functions on $\mathbb{R}$ defined by
$$h_k(x)=(2^k\sqrt{\pi} k!)^{-\frac{1}{2}} H_k(x)e^{-\frac{1}{2}x^2}, \quad k=0, 1, 2, \cdots,$$\vspace{.30mm}
The Hermite functions $\{h_k \}$ are the eigenfunctions of the Hermite operator $H=-\frac{d^2}{dx^2}+x^2$ with eigenvalues $2k+1,  k=0, 1, 2, \cdots$. These functions form an orthonormal basis for $L^2(\mathbb{R})$. The higher dimensional Hermite functions denoted by $\Phi_{\alpha}$ are then obtained by taking tensor product of one dimensional Hermite functions. Thus for any multi-index $\alpha \in \mathbb{N}_0^n$ and $x \in \mathbb{R}^n$, we define
 $\Phi_{\alpha}(x)=\prod_{j=1}^{n}h_{\alpha_j}(x_j).$
The family $\{\Phi_{\alpha}\}$ is then an orthonormal basis for $L^2(\mathbb{R}^n)$. They are eigenfunctions of the Hermite operator $H=-\Delta+|x|^2$ corresponding to eigenvalues $(2|\alpha|+n)$, where $|\alpha |=\sum_{j=1}^{n}\alpha_j$.
\subsection{Pseudo-Differential Operator on $\mathbb{R}^n$}
Given a reasonable function $a$ on $\mathbb{R}^n \times \mathbb{R}^n$, the corresponding operator $T_{a}$ associated  with the function $a $ given by
\begin{align*}
T_{a}f(x)&=a(x, D)f(x)= (2\pi)^{-\frac{n}{2}}\int_{\mathbb{R}^n} e ^{i x\cdot \xi}a(x, \xi)\hat{f}(\xi)\, d\xi,  \quad \forall x \in \mathbb{R}^n
\end{align*} for all Schwartz class functions $f$ on $\mathbb{R}^n$, where the Fourier transform of $f$ is defined by
\begin{align*}
\hat{f}(\xi)&= (2\pi)^{-\frac{n}{2}}\int_{\mathbb{R}^n} {f}(x) e ^{-i x \cdot \xi}\, dx  , \quad \forall \xi \in \mathbb{R}^n.
\end{align*}
The operator $T_a$ is called pseudo-differential operator corresponding to the symbol $a $. Let $m\in\mathbb{R}, 0\leq\delta<\rho\leq 1.$ Then the symbol class $S_{\rho,\delta}^m(\mathbb{R}^n)$ consists of those functions $a(x,\xi)\in C^\infty(\mathbb{R}^n\times\mathbb{R}^n)$ satisfying
\begin{align}
|\partial _{x}^{\alpha }\partial _{\xi}^{\beta } a(x,\xi )| \leq C_{\alpha ,\beta }\,(1+|\xi |^2)^\frac{{m-\delta|\alpha |+\rho|\beta|}}{2}
\end{align}\label{order}
for all multi-indices $\alpha,\beta$. We take $\rho=1$ and $\delta=0$ through out the paper and denote the symbol class by $S^m(\mathbb{R}^n)$.

The Weyl quantization $Op^W$ for a ``reasonable" symbol $a$ in $\mathbb{R}^n \times\mathbb{R}^n$ is given by
$$Op^W(a)f(u)=(2\pi)^{-n} \int_{\mathbb{R}^n}\int_{\mathbb{R}^n} e ^{i (u-v)\cdot \xi}a\left( \frac{u+v}{2},\xi\right){f}(v)\,dv \,  d\xi,  \quad \forall u \in \mathbb{R}^n,$$ for all Schwartz class functions $f$ on $\mathbb{R}^n$. 
The composition of two Weyl quantized operators $Op^W(a)$ and $Op^W(b)$  is  given by $Op^W(a)Op^W(b)=Op^W(a \#   b)$, where  (see \cite{ler})

$$
a \#  b(\zeta, u)=(2 \pi)^{-2 n}  \int_{\mathbb{R}^{n}} \int_{\mathbb{R}^{n}} \int_{\mathbb{R}^{n}} \int_{\mathbb{R}^{n}} e^{-2 i\{(\xi-\zeta)(y-u)-(\eta-\zeta)(x-u)\}} \\
	a(\xi, x) b(\eta, y) d \xi d \eta d x d y
$$
and asymptotically
\begin{align}\label{realerror}
a \#  b (x, \xi)&\sim  \sum_{j=0}^{N} \frac{1}{j !}\left(\frac{i}{2}\right)^{j} a(x, \xi)\left({\overleftarrow{ \frac{\partial }{\partial \xi}} ~ \overrightarrow{\frac{\partial }{\partial u}} -\overrightarrow{\frac{\partial }{\partial \xi}}~ \overleftarrow{\frac{\partial }{\partial u} }}\right)^j b(x, \xi) +S_{N}(x, \xi) 
\end{align}
 (arrows point towards the factor to be differentiated)  with $S_{N} \in S^{m_{1}+m_{2}-N}(\mathbb{R}^n)$.

Further, if $Op^W(a)$ is a trace class operator whose symbol $a(x, \xi) \in L^{1}\left(\mathbb{R}^{n}\times \mathbb{R}^{n}\right)$, then $ {tr} (Op^W(a))=(2\pi)^{-n}\int_{\mathbb{R}^{n}} \int_{\mathbb{R}^{n}} a(x, \xi) d x d \xi$. Moreover,  the correspondence $a \rightarrow Op^W(a) $ is an isometry of $L^{2}\left(\mathbb{R}^{n}\times \mathbb{R}^{n}\right)$ onto the Hilbert-Schmidt operators on $L^{2}\left(\mathbb{R}^{n}\right)$. This yields 
 \begin{align}\label{hilbert1}
 	  {tr} (A B^{*})=\int_{\mathbb{R}^{n}} \int_{\mathbb{R}^{n}} (a \# \bar{b})(x, \xi ) \;d x \;d \xi=\int_{\mathbb{R}^{n}} \int_{\mathbb{R}^{n}}  a(x, \xi) \overline{b(x, \xi)} \;d x \;d \xi,
 \end{align}  
where  $A=Op^W(a)$ and  $B=Op^W(b)$ .

\subsection{Heisenberg Group}
One of the simple and natural example of non-abelian, non-compact group is the famous Heisenberg group $\mathbb{H}^n$, which plays an important role in several branches of mathematics. The Heisenberg group $\mathbb{H}^n$  is a Nilpotent Lie group whose underlying manifold is $ \mathbb{R}^{2n+1} $ and the group operation is defined by
$$(x, y, t)(x', y', t')=(x+x', y+y', t+t'+ \frac{1}{2}(xy'-x'y)),$$
where $(x, y, t)$ and $ (x', y', t')$ are in $\mathbb{R}^n \times \mathbb{R}^n \times \mathbb{R}$. Moreover, $\mathbb{H}^n$ is a unimodular Lie group on which the Haar measure is the usual Lebesgue measure $\, dx \, dy \,dt.$ The canonical basis for the Lie algebra $\mathfrak{h}_n$ of $\mathbb{H}^n$ is given by the left-invariant vector fields:
\begin{align}\label{VF}
X_j&=\partial _{x_j}-\frac{y_j}{2} \partial _{t}, &Y_j=\partial _{y_j}+\frac{x_j}{2} \partial _{t}, \quad j=1, 2, \dotsc   n,   ~ \mathrm{and} ~ T=\partial _{t},
\end{align} \vspace{.30mm}
satisfying the commutator relation
$[X_{j}, Y_{j}]=T, \quad j=1, 2, \dotsc n.$

The sublaplacian and the  full laplacian on the Heisenberg group are  defined as 
$$\mathcal{L}=\sum_{j=1}^{n}(X_j^2+Y_j^2)=\sum_{j=1}^{n}\bigg(\bigg(\partial _{x_j}-\frac{y_j}{2} \partial _{t}\bigg)^2+\bigg(\partial _{y_j}+\frac{x_j}{2} \partial _{t}\bigg)^2\bigg)$$
 and 
 $$\Delta_\mathbb{H}=\sum_{j=1}^{n}(X_j^2+Y_j^2+T_j^2)$$ respectively.
By Stone-von Neumann theorem, the only infinite dimensional unitary irreducible representations
(up to unitary equivalence) are given by $\pi_{\lambda}$, $\lambda$ in $\mathbb{R}^*$, where $\pi_{\lambda}$ is defined by
$$\pi_{\lambda}(x, y, t)f(u)=e^{i \lambda(t+\frac{1}{2}xy)} e^{i \sqrt{\lambda}yu}f(u+\sqrt{|\lambda|x}), \quad \forall u \in \mathbb{R}^n, \forall f  \in L^2(\mathbb{R}^n)\mbox{~and~} (x, y, t)\in \mathbb{H}^n.$$ We use the convention $$\sqrt{\lambda}:={\rm sgn}(\lambda)\sqrt{|\lambda|}
= \begin{cases}  \sqrt{\lambda},&\lambda>0,\\ -\sqrt{|\lambda|}, &\lambda<0.\end{cases}
$$
For each ${\lambda} \in \mathbb{R}^*$, the group Fourier transform of $f\in L^1(\mathbb{H}^n)$ is a bounded linear operator   on  $L^2(\mathbb{R}^n)$ defined by
\begin{align*}
\hat{f}(\lambda) \equiv \pi_{\lambda}(f)=\int_{\mathbb{H}^n}f(x, y, t)\pi_{\lambda}^*(x, y, t)\, dx \, dy \,dt.
\end{align*}We denote  $B(L^2(\mathbb{R}^n))$ to be the set of all bounded operators on $L^2(\mathbb{R}^n)$.
If $f \in L^2(\mathbb{H}^n)$, then $\hat{f}(\lambda)$ is a Hilbert-Schmidt operator on $L^2(\BR^n)$ and satisfies the Plancherel formula $$\int_{\mathbb{R}^*}\|\hat{f}(\lambda)\|_{S_2}^{2} \, d\mu( \lambda)=\|{f}\|_{L^2(\mathbb{H}^n)},$$
where $\| . \|_{S_2}$ stands for  the norm in the Hilbert space $S_2$ of all Hilbert-Schmidt operators  on  $L^2(\mathbb{R}^n)$ and $d\mu(\lambda)=c_n {|\lambda|}^n \, d\lambda$ where $c_n$ is a constant.
\begin{theorem}  For all Schwartz class functions on $\mathbb{H}^n$, the following inversion formula holds:
		\begin{align*}
	f(g) &=  \int_{\mathbb{R}^*}tr~(\pi_{\lambda}(g)\hat{f}(\lambda))\,d\mu(\lambda), \quad \forall g\in \mathbb{H}^n.
	\end{align*}
\end{theorem} 
For a detailed study on the Heisenberg group we refer to Thangavelu \cite{tha98}.
\begin{definition}
	Let $\sigma: \mathbb{H}^n\times \mathbb{R}^* \to B(L^2(\mathbb{R}^n))$  be a operator valued function.  The  pseudo-differential operator $T_\sigma$ corresponding to $\sigma$ is defined by $$	T_\sigma f(g) =  \int_{\mathbb{R}^*}tr~(\pi_{\lambda}(g)\sigma(g, \lambda)\hat{f}(\lambda))\,d\mu(\lambda), \quad  g\in \mathbb{H}^n$$
	for all $f\in \mathcal{S}(\mathbb{H}^n)$.  The operator valued function $\sigma$ is called the symbol of the pseudo-differential operator $T_\sigma$.  We   also often denote  the pseudo-differential operator  $T_\sigma$ as $Op(\sigma)$.
\end{definition}
\section{$(\lambda, V(g))$-Shubin classes $ \Sigma_{\rho, \lambda, V}^{m}\left(\mathbb{R}^{n}\right)$}\label{sec3}
We define the Shubin metric $g_{\xi, u}^{(\rho, \lambda, V(g))}$ depending on  both the parameter $\lambda \in \mathbb{R}^*$ and $V(g), g\in \mathbb{H}^n$ on $\mathbb{R}^{2 n}$  as
$$
g_{\xi, u}^{(\rho, \lambda, V(g))}(d \xi, d u):=\left(\frac{|\lambda|}{1+|\lambda|\left(1+|\xi|^{2}+|u|^{2}\right)+V(g)}\right)^{\rho}\left(d \xi^{2}+d u^{2}\right).
$$
The associated positive function $M^{(\lambda, V(g))}$ on $\mathbb{R}^{2 n}$ is  
$$
M^{(\lambda, V(g))}(\xi, u):=\left(1+|\lambda|\left(1+|\xi|^{2}+|u|^{2}\right)+V(g)\right)^{\frac{1}{2}}.
$$
We consider these $(\lambda, V(g))$-families of metrics for the case $\rho=1 $ as introduced in \cite{BKG}.
\begin{proposition}\label{pro}
	For each $\lambda \in \mathbb{R}^*$ and $g\in \mathbb{H}^n$, the metric $g^{(\rho, \lambda, V(g))}$ is of H\"ormander type  i.e.,  $g$   is uncertain, slowly varying and temperate (see Definition 6.4.2 page 456  of \cite{ruz14}) where the conjugate of $g_{\xi, u}^{(\rho, \lambda, V(g))}$ is $\left(g_{\xi, u}^{(\rho, \lambda, V(g))}\right)^{\omega}$ given by
	$$
	\left(g_{\xi, u}^{(\rho, \lambda, V(g))}\right)^{\omega}(d \xi, d u)=\left(\frac{1+|\lambda|\left(1+|\xi|^{2}+|u|^{2}\right)+V(g)}{|\lambda|}\right)^{\rho}\left(d \xi^{2}+d u^{2}\right).
	$$
	Moreover the gain is  given by
	$$
	\Lambda_{g_{\xi, u}^{(\rho, \lambda, V(g))}}=\left(\frac{1+|\lambda|\left(1+|\xi|^{2}+|u|^{2}\right)+V(g)}{|\lambda|}\right)^{2 \rho}.
	$$
\end{proposition}
\begin{proof}The proof of the proposition follows exactly as in Proposition 1.20 of \cite{BKG} for $\rho=1.$
\end{proof}
For each parameter $\lambda \in \mathbb{R}^*$ and $V(g), g\in \mathbb{H}^n$  we define the $(\lambda, V(g))$-Shubin class  as
$$
\Sigma_{\rho, \lambda, V(g)}^{m}\left(\mathbb{R}^{n}\right):=\{a\in C(\mathbb{R}^n\times\mathbb{R}^n): \|a\|_{\Sigma_{\rho, \lambda, V(g)}^{m}, N}<\infty ~\mbox{for~ each}~ N\in\mathbb{N}_0\},
$$
where
\begin{align*}
&\|a\|_{\Sigma_{\rho, \lambda, V(g)}^{m}, N}\\&\quad :=\sup_{\substack{(\xi, u) \in \mathbb{R}^{n} \times \mathbb{R}^{n}\\|\alpha|,|\beta|\leq N}}|\lambda|^{-\rho \frac{|\alpha|+|\beta|}{2}}\left(1+|\lambda|\left(1+|\xi|^{2}+|u|^{2}\right)+V(g)\right)^{-\frac{m-\rho(|\alpha|+|\beta|)}{2}}\left|\partial_{\xi}^{\alpha} \partial_{u}^{\beta} a(\xi, u)\right|
\end{align*}
is finite. In other words a symbol $a=\{a(\xi, u)\}$ is in $\Sigma_{\rho, \lambda, V(g)}^{m}\left(\mathbb{R}^{n}\right)$ if and only if it satisfies
$$
\begin{array}{c}
\forall \alpha, \beta \in \mathbb{N}_{0}^{n}, \forall(\xi, u) \in \mathbb{R}^{n} \times \mathbb{R}^{n},~~ \mbox{there~exists}~~ C=C_{\alpha, \beta}>0 ~\mbox{such~that}\quad  \\
\left|\partial_{\xi}^{\alpha} \partial_{u}^{\beta} a(\xi, u)\right| \leq C|\lambda|^{\rho \frac{|\alpha|+|\beta|}{2}}\left(1+|\lambda|\left(1+|\xi|^{2}+|u|^{2}\right)+V(g)\right)^{\frac{m-\rho(|\alpha|+|\beta|)}{2}}.
\end{array}
$$

\subsection{The  symbol class $S_{\rho, \delta, \mathcal{H}}^{m}(\mathbb{H}^n)$}\label{subsection1}

We define the  symbol class $S_{\rho, \delta, \mathcal{H}}^{m}(\mathbb{H}^n)$ relative to the operator $1+|\lambda|H+V(g)$  as in Definition 5.2.11 of \cite{ruz14} by the following family of seminorms which are finite:
$$\|\sigma\|_{{S_{\rho, \delta , \mathcal{H}}^{m}}, a, b, c}:=\sup_{g  \in {\mathbb{H}^n}, \lambda \in \mathbb{R}^*} \|\sigma(g, \lambda)\|_{S_{\rho, \delta,, \mathcal{H}, \lambda, V}^{m} a, b, c}, \quad a, b, c \in \mathbb{N}_{0}$$ where
\begin{align}\label{op}
&\|\sigma(g, \lambda)\|_{S_{\rho, \delta,, \mathcal{H}, \lambda, V}^{m} a, b, c}\\\nonumber& \quad  :=\sup_{[\alpha]\leq a, [\beta] \leq b, |\gamma| \leq c}\|(\pi_{\lambda}(I-\mathcal{L})+V(g))^{\frac{\rho[\alpha]-m-\delta[\beta]+\gamma}{2}}X_{g}^{\beta} \Delta'^{\alpha}\sigma(g, \lambda) (\pi_{\lambda}(I-\mathcal{L})+V(g))^{-\frac{\gamma}{2}}\|_{op}.
\end{align}\\
with
$\alpha=(\alpha_1, \alpha_2, \alpha_3)=(\alpha_{11}, \alpha_{12} \dotsc \alpha_{1n}, \alpha_{21}, \alpha_{22} \dotsc \alpha_{2n}, \alpha _3) \in \mathbb{N}_{0}^n \times \mathbb{N}_{0}^n \times \mathbb{N}_{0} , [\alpha] = |\alpha_1|+|\alpha_2|+2\alpha_3$ and $\|\cdot\|_{op}$ denote the operator norm on $B(L^2(\mathbb{R}^n)).$
The difference operators are
$$\Delta'^{\alpha}:=\Delta_{x}^{\alpha_1}\Delta_{y}^{\alpha_2}\Delta_{t}^{\alpha_3}, \quad\mathrm{where}~\Delta_{x}^{\alpha_1}=\Delta_{x_1}^{\alpha_{11}}\Delta_{x_2}^{\alpha_{12}}\cdots \Delta_{x_n}^{\alpha_{1n}},~\Delta_{y}^{\alpha_2}= \Delta_{y_1}^{\alpha_{21}}\Delta_{y_2}^{\alpha_{22}} \dotsc \Delta_{y_n}^{\alpha_{2n}}$$
and
$$X^{\alpha}=X^{\alpha_1}Y^{\alpha_2}T^{\alpha_3},
~\mathrm{where}~ X^{\alpha_1}=X_{1}^{\alpha_{11}}X_{2}^{\alpha_{12}}\dotsc X_{n}^{\alpha_{1n}} ~\mathrm{and} ~ Y^{\alpha_1}=Y_{1}^{\alpha_{11}}Y_{2}^{\alpha_{12}}\cdots Y_{n}^{\alpha_{1n}}.$$

The symbol in  $S_{\rho, \delta, \mathcal{H} }^{m}\left(\mathbb{H}^{n}\right)$ relative to the operator $1+|\lambda|H+V(g)$  can be written  in terms of scalar-valued $(\lambda, V(g))$-symbol. More precisely, the symbols $\sigma = {\sigma(g, \lambda)} $ in $S_{\rho, \delta, \mathcal{H}}^m (\mathbb{H}^n)$ are all of the form $$		\sigma(g, \lambda)=Op^{W}(a_{g, \lambda}(\xi, u)),$$
with the $(\lambda, V(g))$-symbols $a_{g, \lambda}$ 
satisfying some properties described below in terms of the family of $(\lambda, V(g))$-Shubin classes.

\begin{theorem}\label{shubin theorem}
	Let $m, \rho, \delta \in \mathbb{R}$ be real numbers   such that $1\geq \rho\geq \delta\geq 0$ and $(\rho, \delta )\neq(0, 0).$
	
 if $\sigma={\sigma(g, \lambda)}$ is in $S_{\rho, \delta,\mathcal{H}  }^{m}(\mathbb{H}^n)$, then there exist a smooth function $a=a(g, \lambda, \xi, u)=a_{g, \lambda}(\xi, u)$ on $\mathbb{H}^n \times \mathbb{R}^* \times \mathbb{R}^n \times \mathbb{R}^n $ such that
\begin{align*}
	\sigma(g, \lambda)=Op^{W}(a_{g, \lambda})
\end{align*}
	with $\tilde{\partial}_{\lambda, \xi, u}^{\alpha_{3}} X_{g}^{\beta} a_{g, \lambda} \in \Sigma_{\rho, \lambda, V(g)}^{m-2 \rho \alpha_{3}+\delta|\beta|}\left(\mathbb{R}^{n}\right)$ for each $(g, \lambda) \in \mathbb{H}^{n} \times \mathbb{R}^*$ satisfying
\begin{align}\label{shubin estimate}
	\sup _{(g, \lambda) \in \mathbb{H}^{n} \times \mathbb{R}^*  }\left\|\tilde{\partial}_{\lambda, \xi, u}^{\alpha_{3}} X_{g}^{\beta} a_{g, \lambda}\right\|_{\Sigma_{\rho, \lambda, V(g)}^{m-2 \rho \alpha_{3}+\delta|\beta|}\left(\mathbb{R}^{n}\right), N}<\infty
\end{align}
	for every $N \in \mathbb{N}_{0} .$ More precisely, for every $N \in \mathbb{N}_{0}$ there exist $C>0$ and $a, b, c$ such that
$$
	\sup _{(g, \lambda) \in \mathbb{H}^{n} \times \mathbb{R}^*  }\left\|\tilde{\partial}_{\lambda, \xi, u}^{\alpha_{3}} X_{g}^{\beta} a_{g, \lambda}\right\|_{\Sigma_{\rho, \lambda, V(g)}^{m-2 \rho \alpha_{3}+\delta|\beta|}\left(\mathbb{R}^{n}\right), N} \leq C\|\sigma\|_{S_{\rho, \delta, \mathcal{H}}^{m}, a, b, c} .
$$
Conversely, if $a=\{a_{(g, \lambda, \xi, u)}=a_{g, \lambda}(\xi, u)\}$ is a smooth function on $\mathbb{H}^n \times \mathbb{R}^* \times \mathbb{R}^n \times \mathbb{R}^n $ satisfying (\ref{shubin estimate}) for every $N \in \mathbb{N}_{0}$, then there exist a unique symbol $\sigma \in S_{\rho, \delta,\mathcal{H} }^{m}(\mathbb{H}^n)$  such that $
\sigma(g, \lambda)=Op^{W}(a_{g, \lambda}).
$ Furthermore, for every $a, b, c$ there exists $C>0$ and $N \in \mathbb{N}_{0}$ such that
	\begin{align}\label{shubin norm}
	\|\sigma\|_{S_{\rho, \delta, \mathcal{H}}^{m}, a, b, c} \leq C	\sup _{(g, \lambda) \in \mathbb{H}^{n} \times \mathbb{R}^*  }\left\|\tilde{\partial}_{\lambda, \xi, u}^{\alpha_{3}} X_{g}^{\beta} a_{g, \lambda}\right\|_{\Sigma_{\rho, \lambda, V(g)}^{m-2 \rho \alpha_{3}+\delta|\beta|}\left(\mathbb{R}^{n}\right), N} .
		\end{align}
\end{theorem}
\begin{proof}
	The proof  is similar to the proof of  Theorem 6.5.1 of \cite{ruz14}.
\end{proof}
In other words, Theorem \ref{shubin theorem} yields   that   $\sigma\in S_{\rho, \delta,\mathcal{H}  }^{m}(\mathbb{H}^n)$ is equivalent to 	$\sigma(g, \lambda)=Op^{W}(a_{g, \lambda})$ 
  for each $(g, \lambda)\in \mathbb{H}^n \times \mathbb{R}^*$  with  $a_{g, \lambda}\in C^\infty( \mathbb{R}^{2n})$ satisfying: For any $\alpha_1  \in \mathbb{N}_{0}^{2n+1}$ there exists a   constant $C>0$ such that for every $(g, \lambda)\in \mathbb{H}^n \times \mathbb{R}^*$ and for every  $(\xi, u)\in \mathbb{R}^{n}\times \mathbb{R}^{n}$
 	\begin{align*}
 		|\partial_{\xi}^{\alpha}\partial_{u}^{\beta}\tilde{\partial}_{\lambda , \xi, u}^{\bar{\alpha}} X_{g}^{\bar{\beta}}a_{g, \lambda}(\xi, u) | \leq C_{\alpha, \beta, \bar{\alpha}, \bar{\beta}}|\lambda|^{\rho{\frac{|\alpha|+|\beta|}{2}}}(1+|\lambda|(1+|\xi|^2+|u|^2)+V(g))^{\frac{m-\rho|\alpha_1|+\delta|\bar{\beta}|}{2}}.
 	\end{align*}
 We take $\rho=1$ and $\delta=0$ throughout the article and denote the symbol classes
 $ S_{1, 0,\mathcal{H} }^{m}(\mathbb{H}^n)$  by
 $ S_\mathcal{H}^{m}(\mathbb{H}^n)$.
\begin{example}
For  any $\beta \in \mathbb{R}$, $	\pi_{\lambda}(I-\mathcal{L}), V(g)^\beta $ and  $ (1+\lambda^2)^\beta $ are symbols with order $2, 2\beta $ and $4\beta $ respectively.
\end{example}
\begin{remark}\label{ab}
Let $\sigma \in S_\mathcal{H}^{m}(\mathbb{H}^n)$. Then  we have the following properties. \begin{enumerate}
	\item  If $\beta_{o} \in \mathbb{N}_{0}^{n}$ then the symbol $\left\{X_{x}^{\beta_{o}} \sigma(g, \lambda),(g, \lambda)\in \mathbb{H}^n \times \mathbb{R}^* \right\}$ is in $S_\mathcal{H}^{m}(\mathbb{H}^n)$
	  and
	$$
	\left\|X_{g}^{\beta_{o}} \sigma(g, \lambda)\right\|_{S_{\mathcal{H}}^{m }, a, b, c} \leq C_{b, \beta_{o}}\|\sigma(g, \lambda)\|_{S_{\mathcal{H}}^{m}, a, b+\left[\beta_{o}\right], c}.
	$$
	\item If $\alpha_{o} \in \mathbb{N}_{0}^{n}$ then the symbol $\left\{\Delta^{\alpha_{o}} \sigma(g, \lambda),(g, \lambda) \in \mathbb{H}^n \times \mathbb{R}^*  \right\}$ is in $S _\mathcal{H}^{m-\left[\alpha_{o}\right]}(\mathbb{H}^n)$
	  and
	$$
	\left\|\Delta^{\alpha_{o}} \sigma(g, \lambda)\right\|_{S_{\mathcal{H}} ^{m-\left[\alpha_{o}\right]}, a, b, c} \leq C_{a, \alpha_{o}}\|\sigma(g, \lambda)\|_{S_{\mathcal{H}} ^{m}, a+\left[\alpha_{o}\right], b, c}.
	$$
\item If  $ \sigma_1 \in S_\mathcal{H}^{\mu}(\mathbb{H}^n)$ and $\sigma_2 \in   S_\mathcal{H} ^{\nu}(\mathbb{H}^n)$ then $\sigma(g, \lambda)=\sigma_1 (g, \lambda)\sigma_2(g, \lambda)\in S _\mathcal{H}^{\mu+\nu}(\mathbb{H}^n)$ and 
$$\|\sigma(g, \lambda)\|_{{S_{\mathcal{H}}^{\mu+\nu}}, a, b, c}\leq\|\sigma_1(g, \lambda)\|_{{S_{\mathcal{H}}^{\mu}}, a, b, c+a+|\nu|}	\|\sigma_2(g, \lambda)\|_{{S_{\mathcal{H}}^{\nu}}, a, b, c}.$$
    \item If  $ \sigma_1 \in S_\mathcal{H}^{\mu}(\mathbb{H}^n)$ and $\sigma_2 \in   S _\mathcal{H}^{\nu}(\mathbb{H}^n)$ then $\Delta^{\alpha} \sigma_1 X_{x}^{\beta} \sigma_2\in S_\mathcal{H}^{\mu+\nu-[\alpha]}$.
\end{enumerate} 	
\end{remark}

\begin{lemma}\label{CH83}
	If $A$ is a trace class  pseudo-differential operator on $L^2 (\mathbb{H}^n)$  with symbol $a(\cdot, \cdot) \in L^{1}\left(\mathbb{H}^{ n}\times \mathbb{R}^*, S_1, d\mu(\lambda)\right)$,  then
	\begin{align}
		tr(A)=\int_{\mathbb{H}^n } \int_{\mathbb{R}^*} tr(a(g, \lambda)) ~\,d g \,d\mu(\lambda).
	\end{align}
\end{lemma}
\begin{proof}
	For all  $f \in L^2 (\mathbb{H}^n)$, we have
	\begin{align*}
		(Af) (g)&=  \int_{\mathbb{R}^*}{tr}(\pi_{\lambda}^*(g)a(g, \lambda)\hat{f}(\lambda))\,d\mu(\lambda)\\
		&=  \int_{\mathbb{H}^n } \int_{\mathbb{R}^*}{tr}(\pi_{\lambda}^*(g)a(g, \lambda)\pi_{\lambda}(g_1))\,d\mu(\lambda) f(g_1)\,d g_1 \\
		&=  \int_{\mathbb{H}^n } K(g, g_1)f(g_1)\,d g_1,
	\end{align*}
with 
 $$K(g, g_1)=\int_{\mathbb{R}^*}{tr}(\pi_{\lambda}^*(g)a(g, \lambda)\pi_{\lambda}(g_1))\,d\mu(\lambda).$$
Therefore 
	\begin{align*}
		tr(A)= \int_{\mathbb{H}^n } K(g, g) ~\,d g=\int_{\mathbb{H}^n } \int_{\mathbb{R}^*} tr(a(g, \lambda)) ~\,d g \,d\mu(\lambda).
	\end{align*}
\end{proof}
The correspondence $a  \rightarrow Op(a)$ is an isometry from $L^{2}\left(\mathbb{H}^{ n}\times \mathbb{R}^*, S_2, d\mu(\lambda) \right)$ onto the set of  Hilbert-Schmidt operators on $L^2(\mathbb{H}^n)$ via   square integrable kernels \cite{Jiman}. This allows us to write   
\begin{align}\label{schmidt}\nonumber
	{tr} (Op(a )\circ Op(b)^{*})&=\int_{\mathbb{H}^n } \int_{\mathbb{R}^*}  {tr}(a \#_{\mathbb{H}^n} {b}^{(*)})( g, \lambda) \,d g \,d\mu(\lambda)\\
	&=\int_{\mathbb{H}^n } \int_{\mathbb{R}^*}  {tr}(a(g, \lambda) {b}^{(*)}( g, \lambda)) \,d g \,d\mu(\lambda),\end{align}
where $a \#_{\mathbb{H}^n} {b}$ is the symbol of  the composition $Op(a )\circ Op(b)$ (defined in Theorem \ref{hcom})  and   ${b}^{(*)}$ is the symbol of $Op(b)^{*}$, the    adjoint   of $Op(b)$ (see  page 365 of \cite{ruz14}). 

Now as in the proof of Calderón-Vaillancourt theorem (Theorem  5.7.1 of \cite{ruz14}), we  get the following Calderón-Vaillancourt theorem for the symbol class $S_\mathcal{H}^0(\mathbb{H}^n)$. 
\begin{theorem}[The Calderón-Vaillancourt theorem]\label{CH80003}
	Let  
	$\sigma  \in S_\mathcal{H}^{0}(\mathbb{H}^n)$. Then $Op(\sigma )$  extends  a   bounded operator on $L^{2}(\mathbb{H}^n)$.  Moreover, there exist a constant $C>0$ and a seminorm $\|\cdot\|_{S_{\mathcal{ H}} ^{0}, a, b, c}$ with computable integers $a, b, c \in \mathbb{N}_{0}$ independent of $Op(\sigma )$ such that
	$$ \quad\|Op(\sigma ) \phi\|_{L^{2}(\mathbb{H}^n)} \leq C\|\sigma\|_{S _{\mathcal{H}}^{0}, a, b, c}\|\phi\|_{L^{2}(\mathbb{H}^n)}, \quad  \phi \in \mathcal{S}(\mathbb{H}^n) .$$
\end{theorem}
\subsection{Composition of symbols}
 Let  $a\in S_\mathcal{H}^{m_1}(\mathbb{H}^n)$ and $b\in S_\mathcal{H}^{m_2}(\mathbb{H}^n)$. Then the composition of pseudo-differential operators corresponding to the symbols $a$ and $b$ defines a pseudo-differential operator and the symbol $\sigma$ of the composition is given by the following asymptotic expansion (\ref{com}). 
  We add a constraint on $V$ (see \cite{tan}  and \cite{zel83}) which guarantees the asymptotic expansion (\ref{com}).
  \begin{definition}\label{tem}
The potential 	$V$ is said to be temperate potential if there exists $C>0$ such that
	$$ \|(\pi_{\lambda}(I-\mathcal{L})+V(x))^{-1}(\pi_{\lambda}(I-\mathcal{L})+V(xx_1))\|_{op}\leq C|x_1|^k$$
	for all $x,x_1\in \mathbb{H}^n$ and for some constant $k > 0$.
\end{definition}
\begin{theorem}[Composition formula]\label{hcom}
	Let  $a\in S_\mathcal{H}^{m_1}(\mathbb{H}^n)$ and $b\in S_\mathcal{H}^{m_2}(\mathbb{H}^n)$. Then  the composition $Op(a)\circ Op(b)$ is a pseudo-differential operator with symbol $\sigma \in S_\mathcal{H}^{m_1+m_2}(\mathbb{H}^n)$  having asymptotic expansion
	\begin{align}\label{com}
		\sigma(x, \lambda)\sim \sum_{\alpha}\Delta^\alpha a(x, \lambda)X_x^\alpha b(x, \lambda), \vspace{-1.5 cm}
	\end{align}
	where the asymptotic expansion means that for every $M \in \mathbb{N}$, we have
	$$\sigma-\sum_{[\alpha]\leq M}\Delta^\alpha aX_g^\alpha b\in S_\mathcal{H}^{m_1+m_2-M}(\mathbb{H}^n).$$
\end{theorem}
In order to    estimate the reminder term in composition formula, we need  the following lemma.
\begin{lemma}\label{CH811}
 Let $m_{1}, m_{2} \in \mathbb{R}$, $ \beta_{0} \in \mathbb{N}_{0}^{n}$, and $M, M_{1} \in \mathbb{N}_{0} .$ Suppose that 
\begin{align}\label{CH810}
	\begin{cases}
		 m_{2}  \leq 2 M_{1}<M-m_{1} + v_{1} \\
		m_{2} \leq 2 M_{1}<-m_{1}-M.
	\end{cases}
\end{align}
If $M \geq 2 M_{1}$, then only the second condition may be assumed.
Then there exist a constant $C>0$ and two pseudo-norms $\|\cdot\|_{{S_\mathcal{ H}^{m_1,R}}, a_1 , b_1},   \|\cdot\|_{{S_\mathcal{ H}^{m_2}}, 0, b_2, 0} $ such that for any two symbol $a, b$ and for any $(x, \pi) \in \mathbb{H}^n \times \mathbb{R}^*$, we have
\begin{align*}
&	\left\|X_{x}^{\beta_{0}}\left( a \circ b(x, \pi)-\sum_{[\alpha] \leq M} \Delta^{\alpha} a(x, \pi) X_{x}^{\alpha} b(x, \pi)\right)\right\| 
	\leq   C\|a\|_{{S_{\mathcal{H}}^{m_1,R}}, a_1 , b_1}   \|b\|_{{S_{\mathcal{H}}^{m_2}}, 0, b_2, 0}. 
\end{align*}
\end{lemma}
\begin{proof}
	Let $k_1$ and $k_2$ are the kernels of $Op(a)$ and $ Op(b)$ respectively.  	Then we have $Op(a)\circ Op(b)=Op(\sigma),$ where the symbol of the composition is given by $$\sigma(x, \lambda)=\int_{G} k_{1}(x, z) \pi(z)^{*} b \left(x z^{-1}, \lambda\right) d z.$$
	First we  consider the case when $\beta_0=0$.  Thus 
	\begin{align*}
		&	\sigma(x, \lambda)-\sum_{[\alpha] \leq M} \Delta^{\alpha} a (x,  \lambda) X_{x}^{\alpha} b(x,  \lambda) \\
		&	=\int_{\mathbb{H}^n} k_{1}(x, z)  \pi_\lambda(z)^{*}(\pi_{\lambda}(I-\mathcal{L})+V(x))^{M_1}(\pi_{\lambda}(I-\mathcal{L})+V(x))^{-M_1}\\& \quad\quad\quad\quad\quad\times \left(b \left(x z^{-1},  \lambda \right)-\sum_{[\alpha] \leq M} q_{\alpha}\left(z^{-1}\right) X_{x}^{\alpha}  b(x,  \lambda)\right) d z \\
		&	= \sum_{[\beta]=1}^{M_1}\int_{\mathbb{H}^n} k_{1}(x, z)  \pi_\lambda(z)^{*}(\pi_{\lambda}(I-\mathcal{L}))^{\beta}V(x)^{M_1-\beta }(\pi_{\lambda}(I-\mathcal{L})+V(x))^{-M_1}\\& \quad\quad\quad\quad\quad\times \left(b \left(x z^{-1},  \lambda \right)-\sum_{[\alpha] \leq M} q_{\alpha}\left(z^{-1}\right) X_{x}^{\alpha}  b(x,  \lambda)\right) d z \\
		&=\sum_{\left[\beta_{11}\right]+\left[\beta_{22}\right] \leq 2 M_{1}}\sum_{[\beta]=1}^{M_1}\int_{\mathbb{H}^n} \tilde{X}_z^{\beta_{11}}  k_{1}(x, z)V(x)^{M_1-\beta }   \pi_\lambda(z)^{*}  \\& \quad\quad\quad  \times \tilde{X}_z^{\beta_{22}} (\pi_{\lambda}(I-\mathcal{L})+V(x))^{-M_1} R_{x, M}^{b(\cdot,  \lambda)}\left(z^{-1}\right) d z,
	\end{align*}
	where $R_{x, M}^{b(\cdot, \lambda)}(z)=b \left(x z,  \lambda \right)-\sum_{[\alpha] \leq M} q_{\alpha}\left(z\right) X_{x}^{\alpha}  b(x,  \lambda).$
Taking the operator norm on $B(L^2(\mathbb{R}^n))$, we have
	\begin{align*}	&\|	\sigma(x, \lambda)-\sum_{[\alpha] \leq M} \Delta^{\alpha} a (x,  \lambda) X_{x}^{\alpha} b(x,  \lambda)\|_{op}\\
		&\leq \sum_{\left[\beta_{11}\right]+\left[\beta_{22}\right] \leq 2 M_{1}}\sum_{[\beta]=1}^{M_1}\int_{\mathbb{H}^n} \left| \tilde{X}_z^{\beta_{11}}  k_{1}(x, z)V(x)^{M_1-\beta }  \right|     \\& \quad\quad\quad  \times \left\| \tilde{X}_z^{\beta_{22}} (\pi_{\lambda}(I-\mathcal{L})+V(x))^{-M_1} R_{x, M}^{b(\cdot,  \lambda)}\left(z^{-1}\right)\right \|_{op} .
	\end{align*}
	Using  Taylor's estimate for vector-valued functions given in Proposition 3.1.40 and by Theorem 3.1.51 of \cite{ruz14}, there is a constant $c_1$ (depending  on $M$)   such that
	\begin{align*}
		&\left\| \tilde{X}_z^{\beta_{22}} (\pi_{\lambda}(I-\mathcal{L})+V(x))^{-M_1} R_{x, M}^{b(\cdot,  \lambda)}\left(z^{-1}\right)\right \|_{op} \\
		&=	\left\|  (\pi_{\lambda}(I-\mathcal{L})+V(x))^{-M_1} R_{x, M}^{\tilde{X}_z^{\beta_{22}}b(\cdot,  \lambda)}\left(z^{-1}\right)\right \|_{op} \\
		&\leq C_M
		\sum_{\substack{|\gamma| \leq\left(M-\left[\beta_{22}\right]\right) +1\\ |\gamma| > (M-\left[\beta_{22}\right])}}|z|^{[\gamma]} \sup _{|x_{1}|\leq   c_1|z|}\left\|(\pi_{\lambda}(I-\mathcal{L})+V(x))^{-M_1}  X^{\gamma} X^{\beta_{22}}b\left(x x_{1}, \lambda \right)\right\|_{op}\\
		&= C_M
		\sum_{\substack{|\gamma| \leq\left(M-\left[\beta_{22}\right]\right) +1\\ |\gamma| > (M-\left[\beta_{22}\right])}}|z|^{[\gamma]} \sup _{|x_{1}|\leq   c_1|z|}\left\|(\pi_{\lambda}(I-\mathcal{L})+V(x))^{-M_1} (\pi_{\lambda}(I-\mathcal{L})+V(xx_1))^{M_1}\right\|_{op}\\&\quad \quad \quad \times \left\| (\pi_{\lambda}(I-\mathcal{L})+V(xx_1))^{-M_1}  X^{\gamma} X^{\beta_{22}}b\left(x x_{1}, \lambda \right)\right\|_{op}.
	\end{align*}
	
	Using the fact that  $V$ is a temperate potential,  we have
	\begin{align*}
		&\left\| \tilde{X}_z^{\beta_{22}} (\pi_{\lambda}(I-\mathcal{L})+V(x))^{-M_1} R_{x, M}^{b(\cdot,  \lambda)}\left(z^{-1}\right)\right \|_{op} \\
		&\leq C_M c_1^{M_1}
		\sum_{\substack{|\gamma| \leq\left(M-\left[\beta_{22}\right]\right) +1\\ |\gamma| > (M-\left[\beta_{22}\right])}}|z|^{[\gamma]+k M_1} \sup _{|x_{1}|\leq   c_1|z|}\left\| (\pi_{\lambda}(I-\mathcal{L})+V(xx_1))^{-M_1}  X^{\gamma} X^{\beta_{22}}b\left(x x_{1}, \lambda \right)\right\|_{op}.
	\end{align*}
	Let $\sigma_1(x, \lambda)=V(x)^{M_1-\beta } a(x, \lambda).$ Then $\sigma_1\in S_\mathcal{H}^{m_1+2(M_1-\beta)}$ with kernel $\tilde{k}_x=V(x)K_1(x,\cdot)$. So $\sigma_1=\pi(\tilde{k}_x)$.
	Choosing $M, M_1$ such that it satisfies (\ref{CH810}) and the conditions  of Lemma 5.5.6 in  \cite{ruz14}. Thus 
	\begin{align*}	&\|	\sigma(x, \lambda)-\sum_{[\alpha] \leq M} \Delta^{\alpha} a (x,  \lambda) X_{x}^{\alpha} b(x,  \lambda)\|_{op}\\
		&\leq  \sum_{\left[\beta_{11}\right]+\left[\beta_{22}\right] \leq 2 M_{1}}\sum_{\beta=1}^{M_1} C_M c_1^{kM_1}
		\sum_{\substack{|\gamma| \leq\left(M-\left[\beta_{22}\right]\right) +1\\ |\gamma| > (M-\left[\beta_{22}\right])}}\int_{\mathbb{H}^n}|z|^{[\gamma]+k M_1} \left| \tilde{X}_z^{\beta_{11}}  k_{1}(x, z)V(x)^{M_1-\beta }  \right|    d z \\&\quad\quad\quad\quad \times  \|b\|_{{S_{\mathcal{H}}^{m_2}}, 0, b_2, 0}\\
		&\leq  C_2 \|\sigma_1\|_{{S_{\mathcal{H}}^{m_1+2(M_1-\beta ),R_2}}, 0, b_2, 0}  \times  \|b\|_{{S_{\mathcal{H}}^{m_2}}, 0, b_2, 0}\\
		&\leq   C\|a\|_{{S_{\mathcal{H}}^{m_1,R}}, a_1 , b_1}  \times  \|b\|_{{S_{\mathcal{H}}^{m_2}}, 0, b_2, 0},
	\end{align*}
	where $C=C_1\|V\|_{{S_{\mathcal{H}}^{2(M_2-\beta),R}}, a_1 , b_1}$. The  general case $\beta_0\neq 0$   follows  by adopting the   proof of Lemma 5.5.5 in \cite{ruz14}.
\end{proof}
\noindent \textbf{Proof of Theorem \ref{hcom}:}
	Let $T=Op(a)\circ Op(b)$. Then $$Tf(x)=\int_{\mathbb{H}^n}\int_{\mathbb{H}^n}f(z)k_2(y, z^{-1}y)k_1(x, z)dy dz,$$
	where $k_1$ and $k_2$ are the kernels of $Op(a)$ and $ Op(b)$ respectively.   Furthermore, we have $Op(a)\circ Op(b)=Op(\sigma),$ where
	$$\sigma(x, \pi)=\int_{G} k_{1}(x, z) \pi(z)^{*} b \left(x z^{-1}, \pi\right) d z.$$
By the Taylor series expansion (see \cite{ruz14}) of $b$ in the first variable 
	 we get $$\sigma(x, \lambda)\sim \sum_{\alpha}\Delta^\alpha a(x, \lambda)X_x^\alpha b(x, \lambda).$$
The reminder term is estimated  similar to Theorem  5.5.3 of \cite{ruz14} with few modifications. We will only indicate the  main steps with modifications  in our setting. 	Let $m=m_{1}+m_{2}$, 	  $\beta_{0} \in \mathbb{N}_{0}$ and $M_{0} \in \mathbb{N}$. By Theorem \ref{CH80003}, we have
\begin{align}\label{CH8002}\nonumber
	&	\quad \left\|X_{x}^{\beta_{0}}	\tau_{M}(x, \pi)(\pi_{\lambda}(I-\mathcal{L})+V(x))^{-\frac{m-  M_{0} }{2}}\right \|_{op}\\\nonumber
	&=\left 	\|X_{x}^{\beta_{0}}	\tau_{M}(x, \pi)(\pi_{\lambda}(I-\mathcal{L}))^{-\frac{m-  M_{0} }{2}}\left[ (\pi_{\lambda}(I-\mathcal{L}))^{\frac{m-  M_{0} }{2}} (\pi_{\lambda}(I-\mathcal{L})+V(x))^{-\frac{m-  M_{0} }{2}}\right]\right \|_{op}\\\nonumber
	&\leq \left 	\|X_{x}^{\beta_{0}}	\tau_{M}(x, \pi)(\pi_{\lambda}(I-\mathcal{L}))^{-\frac{m-  M_{0} }{2}}\right\|\left\|(\pi_{\lambda}(I-\mathcal{L}))^{\frac{m-  M_{0} }{2}} (\pi_{\lambda}(I-\mathcal{L})+V(x))^{-\frac{m-  M_{0} }{2}}\right \|_{op}\\
	&\leq C\left 	\|X_{x}^{\beta_{0}}	\tau_{M}(x, \pi)(\pi_{\lambda}(I-\mathcal{L}))^{-\frac{m-  M_{0} }{2}}\right\|_{op},
\end{align}
where $\tau_{M}=a\circ b-\sum_{[\alpha]\leq M}\Delta^\alpha a X_x^\alpha b $.
%
We fix $m_{2}^{\prime}:=-m_{1}+  M_{0}$. Then we can find $M \geq \max \left(M_{0}, v_{1}\right)$ such that $ 	-m_{1}+M- m_{2}^{\prime}  \geq 2.$	This shows that we can find $M_{1}$ satisfying the second condition of  
(\ref{CH810})  for $m_{1}, m_{2}^{\prime}$ and therefore also the first. Hence we can apply Lemma  
\ref{CH811}
 to $M, M_{1}$ and the symbols $a$ and $b(\pi_{\lambda}(I-\mathcal{L}))^{-\frac{m-M_{0} }{2}}$, with orders $m_{1}$ and $m_{2}^{\prime} .$  
Thus by (\ref{CH8002})   and   Theorem \ref{CH80003}, we have 
\begin{align*}
	 & \left\|X_{x}^{\beta_{0}}	\tau_{M}(x, \pi)(\pi_{\lambda}(I-\mathcal{L})+V(x))^{-\frac{m-  M_{0} }{2}}\right \|_{op}\\
	\leq& 	\left\|a\right\|_{S_{\mathcal{H}}^{m_{1}, R}, a_{1}, b_{1}}\left  \| b~(\pi_{\lambda}(I-\mathcal{L}))^{-\frac{m-  M_{0} }{2}}\right \|_{S_{\mathcal{H}}^{m_2'}, 0, b_{2}, 0}\\
	=& 	\left\|a\right\|_{S_{\mathcal{H}}^{m_{1}, R}, a_{1}, b_{1}} \left \|  (\pi_{\lambda}(I-\mathcal{L}+V(g)))^{-\frac{m_2' }{2}}  b(\pi_{\lambda}(I-\mathcal{L}))^{-\frac{m-  M_{0} }{2}}\right  \|_{op}\\
	=& 	\left\|a\right\|_{S_{\mathcal{H}}^{m_{1}, R}, a_{1}, b_{1}} \left\|(\pi_{\lambda}(I-\mathcal{L}+V(g)))^{-\frac{m_2' }{2}}  b~ (\pi_{\lambda}(I-\mathcal{L}+V(g)))^{-\frac{m-  M_{0} }{2}}\right. \\&\left. \times  (\pi_{\lambda}(I-\mathcal{L}+V(g)))^{-\frac{m-  M_{0} }{2}} (\pi_{\lambda}(I-\mathcal{L}))^{-\frac{m-  M_{0} }{2}} \right\|_{op}\\
		\leq & 	\left\|a\right\|_{S_{\mathcal{H}}^{m_{1}, R}, a_{1}, b_{1}}\left  \|  (\pi_{\lambda}(I-\mathcal{L}+V(g)))^{-\frac{m_2' }{2}}  b~ (\pi_{\lambda}(I-\mathcal{L}+V(g)))^{-\frac{m-  M_{0} }{2}}\right \|_{op}\\&\times \left\| (\pi_{\lambda}(I-\mathcal{L}+V(g)))^{-\frac{m-  M_{0} }{2}} (\pi_{\lambda}(I-\mathcal{L}))^{-\frac{m-  M_{0} }{2}}\right  \|_{op}\\
			\leq & C	\left\|a\right\|_{S_{\mathcal{H}}^{m_{1}, R}, a_{1}, b_{1}}\left  \|    b~ (\pi_{\lambda}(I-\mathcal{L}+V(g)))^{-\frac{m-  M_{0} }{2}}\right \|_{S_{\mathcal{H}}^{m_2'}, 0, b_{2}, 0}\\
				\leq & C	\left\|a\right\|_{S_{\mathcal{H}}^{m_{1}, R}, a_{1}, b_{1}}\left  \|    b \right \|_{S_{\mathcal{H}}^{m_2}, 0, b_{2}, c_2}.
\end{align*}

The rest of proof follows along the  similar lines of  Theorem  5.5.3 in \cite{ruz14}.

\section{Symbolic calculus relative to $(1+|\lambda|H+V(g)+w)$ on the Heisenberg group}\label{sec4}
Let $\Gamma$   $ \subset \mathbb{C}$ be a curve enclosing $\mathbb{R}^{+}$ and $w$ vary over  $\Gamma$.  In particular, let us consider  the curve $\Gamma$ be made up of two half-lines hinged at $-1$ and makes  an angles of $\pm \frac{\pi}{4}$ with respect to  the real axis. In order to construct the pseudo-differential approximation to the operator $(\mathcal{ H}+u)^{-m}$, we need to define the following symbol class.
\subsection{The  symbol class $S_{\rho, \delta,\mathcal{H}, w}^{m}(\mathbb{H}^n)$}
We define the  symbol class $S_{\rho, \delta,\mathcal{H}, w}^{m}(\mathbb{H}^n)$ relative to the operator $1+|\lambda|H+V(g)+|w|$    defined as  in Subsection \ref{subsection1}  and $(\lambda, V(g))$-Shubin class defined as  in Section \ref{sec3} relative to the weight $1+|\lambda|(|\xi|^2+|x|^2+1)+V(g)+|w|$.  Also we get the similar  result for the symbol class $S_{\rho, \delta,\mathcal{H},w}^{m}(\mathbb{H}^n)$ as in Theorem \ref{shubin theorem}.  When  $\rho=1$ and $\delta=0$, we denote   the symbol classes $S_{1, 0,\mathcal{H}, w}^{m}(\mathbb{H}^n)$ by $S_{\mathcal{H}, w}^{m}(\mathbb{H}^n)$.

\begin{proposition}
	Let $a_{g, \lambda, w}(\xi, u)=\left(|\lambda|(1+|\xi|^2+|u|^2)+V(g)-w\right)^{s}, s\in \mathbb{R} $ and   $\sigma(g, \lambda, w)$ $=Op^{W}(a_{g, \lambda, w}).$ Then $\sigma\in S_{\mathcal{H}, w}^{2s}(\mathbb{H}^n)$.
\end{proposition}
\begin{proof}
By Theorem \ref{shubin theorem}, 	   $Op^{W}( |\lambda|(1+|\xi|^2+|u|^2)+V(g)-w ) \in S_{\mathcal{H},w}^{2}(\mathbb{H}^n)$. Now 
	 \begin{align*}
	 &	\partial_{\xi}^{\alpha} \partial_{u}^{\beta} \tilde{\partial}_{\lambda, \xi, u}^{\tilde{\alpha}} X_{g}^{\tilde{\beta}} a_{g, \lambda, w}(\xi, u) \\&\quad\quad =\sum_{\substack{1\leq \theta \leq |\alpha|+|\beta|+|\tilde{\alpha}|+|\tilde{\beta}|\\|\mu_1|+\cdots+|\mu_a|\leq |\alpha|\\|\nu_1|+\cdots+|\nu_a|\leq |\beta|\\|\tilde{\mu_1}|+\cdots+|\tilde{\mu_a}|\leq |\tilde{\alpha}|\\|\tilde{\nu_1}|+\cdots+|\tilde{\nu_a}|\leq |\tilde{\beta}|}} \left(|\lambda|(1+|\xi|^2+|u|^2)+V(g)-w\right)^{s-\theta}\\&\quad\quad \times  \prod_{j=1}^{\theta}\partial_{\xi}^{\mu_j} \partial_{u}^{\nu_j} \tilde{\partial}_{\lambda, \xi, u}^{\tilde{\mu_j}} X_{g}^{\tilde{\nu_j}} \left(|\lambda|(1+|\xi|^2+|u|^2)+V(g)-w\right).
	 \end{align*}
	 Since each term is bounded by a constant times  \begin{align*}
&	 \left(|\lambda|(1+|\xi|^2+|u|^2)+V(g)-w\right)^{s-\theta}   \prod_{j=1}^{\theta }|\lambda|^{\frac{|\mu_j|+|\nu_j|}{2}}\\&\quad \quad \times (1+|\lambda|(1+|\xi|^2+|u|^2)+V(g)+|w|)^{\frac{2-2 |\tilde{\mu_j}| - (|\mu_j|+|\nu_j|)}{2}}\\
&\leq |\lambda|^{\frac{|\alpha|+|\beta|}{2}}(1+|\lambda|(1+|\xi|^2+|u|^2)+V(g)+|w|)^{\frac{-2-2 |\tilde{\alpha}| - (|\alpha|+|\beta|)}{2}},
  \end{align*}
	 thus  for  any  $w\in \Gamma$, we have 
	  \begin{align*}
	 	\left| 	\partial_{\xi}^{\alpha} \partial_{u}^{\beta} \tilde{\partial}_{\lambda, \xi, u}^{\tilde{\alpha}} X_{g}^{\tilde{\beta}} a_{g, \lambda, w}(\xi, u) \right| \leq C |\lambda|^{\frac{|\alpha|+|\beta|}{2}}(1+|\lambda|(1+|\xi|^2+|u|^2)+V(g)+|w|)^{\frac{-2-2 |\tilde{\alpha}| - (|\alpha|+|\beta|)}{2}}.
	 \end{align*}
	 Now 
	 \begin{align*}
	 &	\left\|\tilde{\partial}_{\lambda, \xi, u}^{\tilde{\alpha} }X_{g}^{\tilde{\beta}} a_{g, \lambda,w}\right\|_{\Sigma_{\rho, \lambda}^{-2-2   \tilde{\alpha}}\left(\mathbb{R}^{n}\right), N}\\
	 &=\sup_{\substack{(\xi, u) \in \mathbb{R}^{n} \times \mathbb{R}^{n}\\|\alpha|,|\beta|\leq N}}|\lambda|^{- \frac{|\alpha|+|\beta|}{2}}\left(1+|\lambda|\left(1+|\xi|^{2}+|u|^{2}\right)+V(g)+w\right)^{-\frac{-2-2|\tilde{\alpha}|- (|\alpha|+|\beta|)}{2}}\\&\qquad \times\left|\partial_{\xi}^{\alpha} \partial_{u}^{\beta}\tilde{\partial}_{\lambda, \xi, u}^{\tilde{\alpha} }X_{g}^{\tilde{\beta}} a_{g, \lambda,w}(\xi, u)\right|	\leq C_{\tilde{\alpha}, \tilde{\beta}, N}.
	 \end{align*}
 Thus	 $\sigma=\sigma(g, \lambda, w)=Op^{W}(a_{g, \lambda, w}) \in S_{ \mathcal{H}, w}^{2s}(\mathbb{H}^n)$ by (\ref{shubin norm}).
\end{proof}
Construct  a symbol $R_N(g, \lambda, w)$ such  that $(\mathcal{H}-w)\circ Op(R_N(g, \lambda, w)) =I_{L^2(\mathbb{H}^n)}+Op(S_N(g, \lambda, w))$, where $S_N \in S_{\mathcal{H},w}^{-N}(\mathbb{H}^n)$ or equivalently  
$ (|\lambda|(H+I)+V(g)-w)\#_{\mathbb{H}^n}R_N(g, \lambda, w)$ $=I_{L^2(\mathbb{R}^n)}+S_N(g, \lambda, w)$.
 By substituting  the expansion $R_N = R_{-2}+R_{-3} + \cdots  + R_{-N}$ with the property that $R_{-2-\ell } \in S_{\mathcal{H},w}^{-2-\ell}(\mathbb{H}^n)$  into the asymptotic expansion (\ref{com}), we get 
\begin{align*}
&(|\lambda|(H+I)+V(g)-w)	\#_{\mathbb{H}^n} R_N(g, \lambda, w) \\&=\sum_{[\alpha]\leq N}\Delta^\alpha  (|\lambda|(H+I)+V(g)-w)	 ~X_g^\alpha R_N(g, \lambda, w)+S_N(g, \lambda, w)
	\\&=I_{L^2(\mathbb{R}^n)}+S_N(g, \lambda, w).
\end{align*}
Now  solving  for $R_{-2-\ell}$ recursively by comparing the  order by  order of the symbols so that the sum equals to 1, we get 
\begin{align}\label{CH81}
	R_{-2} (g, \lambda, w)=(|\lambda|(H+I)+V(g)-w)^{-1}
\end{align} and 
\begin{align}\label{R2}\nonumber
&	R_{-2-\ell} (g, \lambda, w)\\& = (|\lambda|(H+I)+V(g)-w)^{-1} \sum_{\substack{|k|+|j|= \ell \\|k|<\ell}}\Delta^j  (|\lambda|(H+I)+V(g)-w)	~ X_g^j R_{-2-|k|}(g, \lambda, w) 
\end{align}
for $w\in \Gamma.$ To understand the dependence on $r$, we express   the symbol   $R_{-2-\ell} $   differently in the following proposition.
\begin{proposition}
 Let $w\in \Gamma$. Then  
\begin{align}\label{CH82} 
R_{-2-\ell} (g, \lambda, w)&=(|\lambda|(H+I)+V(g)-w)^{-1}\sum_{[\frac{\ell}{2}] \leq M\leq \ell }  R_{\ell, M}(g, \lambda)(|\lambda|(H+I)+V(g)-w)^{-M} ,
\end{align}
  where $[\frac{\ell}{2}]$ denotes  the least integer   grater than  $\frac{\ell}{2}$  and    $R_{\ell, M}(g, \lambda) \in   S_{\mathcal{H},w}^{2 M-\ell}(\mathbb{H}^n)$ is  a polynomial in  $\pi(X) $ and $X_g^\alpha V ,|\alpha| \leq \ell$.
\end{proposition}\begin{proof}
 We prove the proposition by induction on $\ell$. When  $\ell=0$, the expression is trivial from (\ref{CH81}). Assume that the expression (\ref{CH82}) holds for $k \leq \ell-1$.  From (\ref{R2}), the  difference  operator $\Delta$   contributes  only some possible factors of $\pi(X)$ but  no $w$. However,  
 the differential operator $X_g$  either   acts on  $(|\lambda|(H+I)+V(g)-w)^{-M}$ or   $R_{\ell, M}$ (after substituing (\ref{CH82}) for $k\leq \ell-1$ in (\ref{R2}))  resulting  the   expressions as in (\ref{CH82}). It is easy to check that each term in the sum for $R_{-2-\ell}$ lies in $S_{\mathcal{H},w}^{-2-\ell}(\mathbb{H}^n)$ after expanding    by  Leibniz rule.      Since    $ (|\lambda|(H+I)+V(g)-w)^{-1}  R_{\ell, M}(g, \lambda)(|\lambda|(H+I)+V(g)-w)^{-M} \in S_{\mathcal{H},w}^{-2-\ell}(\mathbb{H}^n)$,    $R_{\ell, M} \in S_{\mathcal{H},w}^{2M-\ell}(\mathbb{H}^n)$.  So $R_{\ell, M}$  is     a   polynomial in  $\pi(X) $ and $X_g^\alpha V.$ The highest power of $(|\lambda|(H+I)+V(g)-w)^{-1} $ in the right comes  out  when we throw all derivatives on factors of $(|\lambda|(H+I)+V(g)-w)^{-1} $ and count this number which is essentially  $\ell$.
 
\end{proof}
\subsection{Approximation of symbols}\label{subsec4.1}
Let $f$ be a  holomorphic  function. Then by   the holomorphic functional calculus for unbounded operators, all pseudo-differential approximations can be written in the  following way:
\begin{align*}
	f(\mathcal{H})&=\frac{1}{2\pi i}\int_{\Gamma} f(w)(\mathcal{H}-w)^{-1} \;d w
\end{align*}
and define the pseudo-differential operator $ f_{N}(\mathcal{H})=\frac{1}{2\pi i}\int_{\Gamma} f(w)Op(R_N(g, \lambda, w)) \;d w$ with symbol
$$
\begin{aligned}
	f_{N}(g, \lambda) &=\sum_{\ell =0}^{N} \sum_{[\frac{\ell}{2}] \leq M\leq \ell }  \bigg(   \frac{1}{2\pi i} \int_{\Gamma} f(w)(|\lambda|(H+I)+V(g)-w)^{-1} R_{\ell, M}(g, \lambda )\bigg.\\
	&\bigg.\qquad \times(|\lambda|(H+I)+V(g)-w)^{-M} \;dw \bigg) 
\end{aligned}
$$
by formally computing the residue.  However, the error term  in the  approximation of  $f(\mathcal{H})$ is given by  
\begin{align}\label{errorterm}
	\frac{1}{2\pi i}\int_{\Gamma} f(w)(\mathcal{H}-w)^{-1} Op(S_N(g, \lambda,w))  \;d w.
\end{align}
In particular, letting  $f(w)=(w+u)^{s}$ for some fixed $u>0$,    the   pseudo-differential  approximation to  $(\mathcal{H}+u)^s$  is $Op(	\sigma_{s,N}(g, \lambda))$, where 
$$ \begin{aligned}
	\sigma_{s,N}(g, \lambda) =&\sum_{\ell =0}^{N} \sum_{[\frac{\ell}{2}] \leq M\leq \ell }  \bigg(   \frac{1}{2\pi i} \int_{\Gamma} (w+u)^s(|\lambda|(H+I)+V(g)-w)^{-1} R_{\ell, M}(g, \lambda )\bigg.\\
	&\bigg.\qquad \times(|\lambda|(H+I)+V(g)-w)^{-M} \;dw \bigg).
\end{aligned}
$$
Let $A$ be a linear operator on $L^2(\mathbb{H}^n)$.  Observe that $|tr(A)|\leq  \|(I+\mathcal{ H})^{s}A\| |tr(I+\mathcal{ H})^{-s}|$. If $\|(I+\mathcal{ H})^{s}A\|$ is finite and $(I+\mathcal{ H})^{-s}$ is a trace class operator the $A$ is a trace class operator on $L^2(\mathbb{H}^n)$. Taking $f(w)=(w+1)^s$ in the above discussion,  we get $(I+\mathcal{H})^{-\frac{s}{2}}=Op((1+|\lambda|H+V(g))^{-\frac{s}{2}} )+Op(F_{\frac{s}{2}}(g, \lambda))$. So  $(I+\mathcal{ H})^{-s}$ is a  trace class operator when $(I+\mathcal{ H})^{-\frac{s}{2}}$ is Hilbert-Schmidt operator.  That means,  if   $Op((1+|\lambda|H+V(g))^{-\frac{s}{2}} )$ and $Op(F_{\frac{s}{2}}(g, \lambda))$ are   Hilbert-Schmidt operators or 
equivalently  $(1+|\lambda|H+V(g))^{-\frac{s}{2}} , F_{\frac{s}{2}}(g, \lambda) \in L^{2}\left(\mathbb{H}^{ n}\times \mathbb{R}^*, S_2, d\mu(\lambda) \right)$, $(I+\mathcal{H})^s$ is a trace class operator.   
By  generalized Minkowski's   inequality, we have 
	\begin{align}\label{CH1000}\nonumber
		&	\|(1+|\lambda|H+V(g))^{-\frac{s}{2}}\|_{L^{2}\left(\mathbb{H}^{ n}\times \mathbb{R}^*, S_2, d\mu(\lambda) \right)}\\\nonumber
		&=\left(\int_{\mathbb{H}^n } \int_{\mathbb{R}^*} |tr(1+|\lambda|H+V(g))^{-\frac{s}{2}}|^2 \,d g \,d\mu(\lambda) \right)^{\frac{1}{2}}\\\nonumber
		&=\left(\int_{\mathbb{H}^n } \int_{\mathbb{R}^*} \left| \sum_{\alpha} \frac{1}{(1+|\lambda| (2|\alpha|+n)+V(g))^{\frac{s}{2}}}\right |^2 \,d g \,d\mu(\lambda) \right)^{\frac{1}{2}}\\\nonumber
		&\leq C\sum_{\alpha}  \left(\int_{\mathbb{H}^n } \int_{\mathbb{R}^*}  \frac{\lambda^n}{(1+|\lambda| (2|\alpha|+n)+V(g))^{s}} \,d g \,d\lambda \right)^{\frac{1}{2}}\\\nonumber
		&=C\sum_{\alpha}  \frac{1}{(2|\alpha|+n)^{\frac{n+1}{2}}  }  \left(\int_{\mathbb{H}^n } 	\int_{1+V(g)}^\infty    \frac{(u-1-V(g))^n}{u^{s}} \, du\,d g   \right)^{\frac{1}{2}}\\
		&=C\sum_{\alpha}  \frac{1}{(2|\alpha|+n)^{\frac{n+1}{2}}  }  \left(\int_{\mathbb{H}^n } \frac{1}{(1+V(g))^{s-n-1}} \,d g   \right)^{\frac{1}{2}}.
	\end{align}
	Under the assumption $V(g) \sim V_{0}|g|^{k}$ as $|g| \rightarrow \infty$,  the function $\frac{1}{(1+V(g))^{s-n-1}}$  is integrable if we choose   $(s-n-1) k>1$. A   similar   argument gives  $F_{\frac{s}{2}}$ is also a Hilbert-Schmidt  operator for large $N$.  Indeed $(\mathcal{H}-w)^{-1}=Op(R_N)+(\mathcal{H}-w)^{-1}Op(S_N(g, \lambda, w))$  is compact and hence has discrete spectrum.
 
\begin{proposition}\label{CH86}
	Let $u>0$ and  $m\in \mathbb{N}$ be  sufficiently large such that  $	(\mathcal{H}+u)^{-m}$ is  a trace class operator on $L^2(\mathbb{H}^n)$. Then  for such $m$, 
$	(\mathcal{H}+u)^{-m}=Op\left( (|\lambda|(H+I)+V(g)+u)^{-m}\right)+Op(E(g, \lambda, u))$ 
such that 
\begin{align*}
&	\left| 	tr(\mathcal{H}+u)^{-m}-tr \left( Op\left( (|\lambda|(H+I)+V(g)+u)^{-m}\right)\right) \right|= \left|tr(Op(E(g, \lambda, u)))\right| \\&\leq  \psi_1(u)\left|tr \left( Op\left( (|\lambda|(H+I)+V(g)+u)^{-m}\right)\right)\right|
\end{align*}
with 
 $\psi_1(u)\to 0  $ as $u\to \infty.$
\end{proposition}

\begin{proof}
From the discussions in the previous  subsections,  we   write  
 \begin{align}\label{CH812}
 	(\mathcal{H}+u)^{-m}=Op\left( (|\lambda|(H+I)+V(g)+u)^{-m}\right)+Op(E(g, \lambda , u)),
 \end{align}
 where
	\begin{align*}
		&E(g, \lambda, u)\\&=\sum_{\ell =1}^{N} \sum_{[\frac{\ell}{2}] \leq M\leq \ell }         \frac{1}{2\pi i} \int_{\Gamma} (w+u)^{-m}(|\lambda|(H+I)+V(g)-w)^{-1} R_{\ell, M}(g, \lambda )  \\
		&  \times(|\lambda|(H+I)+V(g)-w)^{-M} \;dw     +		\frac{1}{2\pi i}\int_{\Gamma} (w+u)^{-m}(\mathcal{H}-w)^{-1} Op(S_N(g, \lambda, w))  \;d w.
	\end{align*}
  For large $N$, choose   $0<s<N$ such that  $(I+\mathcal{ H})^{-\frac{s}{2}}$ is a trace class operator. Then    $Op(S_N(g, \lambda, w))$ is a trace class operator with 
\begin{align}\label{CH87}\nonumber
	\left| {tr} (Op(S_N(g, \lambda, w)))\right|
	 &\leq	\left| {tr} \left((I+\mathcal{ H})^{-\frac{s}{2}}(I+\mathcal{ H})^{\frac{s}{2}} Op(S_N(g, \lambda, w))\right)\right| \\
	 &\leq	\left| {tr} \left((I+\mathcal{ H})^{-\frac{s}{2}}\right)\right| \left\|  (I+\mathcal{ H})^{\frac{s}{2}} Op(S_N(g, \lambda, w))\right\|.
\end{align}
But from  Theorem \ref{CH80003}, we have 
\begin{align}\label{CH88}\nonumber
&  \left\|  (I+\mathcal{ H})^{\frac{s}{2}} Op(S_N(g, \lambda, w))\right\|\\\nonumber
&=  \left\|  (I+\mathcal{ H})^{\frac{s}{2}}   (\mathcal{ H}+|w|)^{-\frac{N}{2}}   (\mathcal{ H}+|w|)^{\frac{N}{2}} Op(S_N(g, \lambda, w))\right\|\\\nonumber
  &=  \left\|  (I+\mathcal{ H})^{\frac{s}{2}}(\mathcal{ H}+|w|)^{-\frac{N}{2}}   \right\|  \left\|  (\mathcal{ H}+|w|)^{\frac{N}{2}} Op(S_N(g, \lambda, w))\right\| \\\nonumber
    &\leq C \left\|  (I+\mathcal{ H})^{\frac{s}{2}} (\mathcal{ H}+|w|)^{-\frac{s}{2}} (\mathcal{ H}+|w|)^{\frac{s}{2}} (\mathcal{ H}+|w|)^{-\frac{N}{2}}  \right\|  \\\nonumber
        &\leq C \left\|  (I+\mathcal{ H})^{\frac{s}{2}} (\mathcal{ H}+|w|)^{-\frac{s}{2}} \right\|   \left\|  (\mathcal{ H}+|w|)^{\frac{(s-N)}{2}}   \right\| \\ 
  &=O\left(|w|^{\frac{(s-N)} {2}}\right).
\end{align}
Therefore, from (\ref{CH87}) and (\ref{CH88}), we obtain
\begin{align*}
	\left| {tr}\left( \int_{\Gamma} (w+u)^{-m}(\mathcal{H}-w)^{-1} Op(S_N(g, \lambda, w))  \;d w\right) \right|  \leq C u^{1-m}\to 0~\text{as}~ u\to \infty.
\end{align*}
Consequently this part of the error is negligible and the pseudo-differential part of $E(g, \lambda, u)$ is a trace class operator because it has smooth rapidly decaying symbol.  By Lemma \ref{CH83}, each term of $tr(Op(E(g, \lambda, u))$  is of the form    
\begin{align*}
	& tr\bigg[  Op\bigg( \frac{1}{2\pi i} \int_{\Gamma} (w+u)^{-m}(|\lambda|(H+I)+V(g)-w)^{-1} R_{\ell, M}(g, \lambda )\bigg.\\\bigg.&\qquad\times  (|\lambda|(H+I)+V(g)-w)^{-M} dw \bigg)\bigg]\\
&=\int_{\mathbb{H}^n } \int_{\mathbb{R}^*}	 {tr}\bigg[ \left( \frac{1}{2\pi i} \int_{\Gamma} (w+u)^{-m}(|\lambda|(H+I)+V(g)-w)^{-1-M}  dw\right) R_{\ell, M}(g, \lambda )  \bigg]\,d g \,d\mu(\lambda)\\
&=C_{m, M}\int_{\mathbb{H}^n } \int_{\mathbb{R}^*}  {tr}\left(  { (|\lambda|(H+I)+V(g)+u)^{-m-M}}R_{\ell, M}(g, \lambda )\right) \,d g \,d\mu(\lambda)\\
&=C_{m, M}u^{-m-M}\int_{\mathbb{H}^n } \int_{\mathbb{R}^*}  {tr}\left( { (u^{-1} |\lambda|(H+I)+u^{-1}V(g)+1)^{-m-M}}R_{\ell, M}(g, \lambda ) \right)\,d g \,d\mu(\lambda)\\
&\sim C_{m, M} u^{-m-M+n+1+\frac{2n}{\kappa}+\frac{2}{\kappa}}\int_{\mathbb{H}^n } \int_{\mathbb{R}^*}  {tr}\left(  { ( |\lambda|(H+I)+|g|^\kappa+1)^{-m-M}}R_{\ell, M}(\tilde{g}, u\lambda )\right)|\lambda|^n\,d g \,d \lambda,
\end{align*}
where $ \widetilde{g}=(u^\frac{1}{\kappa}x_1, u^\frac{1}{\kappa}x_2,\cdots, u^\frac{1}{\kappa}x_{2n},u^\frac{2}{\kappa}t)$.  Since $  R_{\ell, M}\in S_{\mathcal{ H}, w}^{2M-\ell }(\mathbb{H}^n)$,  $$\left\| R_{\ell, M}(\tilde{g}, u\lambda ) ( u|\lambda|(H+I)+u|\tilde{g}|^k+1)^{\frac{-2M+\ell}{2}} \right\|_{op}$$ is uniformly bounded and so
\begin{align}\label{approximation}\nonumber
	&\left| 	 {tr}\left( Op\left(  { (|\lambda|(H+I)+V(g)+u)^{-m-M}}R_{\ell, M}(g, \lambda ) \right) \right)\right|\\\nonumber
&\sim  C\frac{u^{n+1+\frac{2n}{\kappa}+\frac{2}{\kappa}}}{u^{m+M}} \int_{\mathbb{H}^n } \int_{\mathbb{R}^*} \bigg|  {tr}\bigg(      ( u|\lambda|(H+I)+u|g|^\kappa+1)^{\frac{2M-\ell}{2}}\bigg.\\&\bigg.	\qquad \times { ( |\lambda|(H+I)+|g|^\kappa+1)^{-m-M}}	 \bigg) \bigg| |\lambda|^n d g \,d \lambda\\\nonumber
&\leq C u^{-m-M+n+1+\frac{2n}{\kappa}+\frac{2}{\kappa}}\int_{\mathbb{H}^n } \int_{\mathbb{R}^*} \sum_{k} \frac{( u|\lambda|(2|k|+n+1)+u|g|^\kappa+1)^{\frac{2M-\ell}{2}}}	{ ( |\lambda|(2|k|+n+1)+|g|^\kappa+1)^{m+M}}|\lambda|^n\	  \,d g \,d \lambda\\\nonumber
&\leq C u^{-m-M+n+1+\frac{2n}{\kappa}+\frac{2}{\kappa}+\frac{2M-\ell}{2}}\sum_{k}\int_{\mathbb{H}^n } \int_{\mathbb{R}^*}  \frac{( |\lambda|(2|k|+n+1)+|g|^\kappa+1)^{\frac{2M-\ell}{2}}}	{ ( |\lambda|(2|k|+n+1)+|g|^\kappa+1)^{m+M}}	|\lambda|^n\  \,d g \,d \lambda\\\nonumber
&\leq C u^{-m+n+1+\frac{2n}{\kappa}+\frac{2}{\kappa}-\frac{\ell}{2}}\sum_{k} \int_{\mathbb{H}^n } \int_{\mathbb{R}^*}  \frac{ 	|\lambda|^n}	{ ( |\lambda|(2|k|+n+1)+|g|^\kappa+1)^{m+\frac{\ell}{2}}}\  \,d g \,d \lambda\\\nonumber
&\leq C u^{-m+n+1+\frac{2n}{\kappa}+\frac{2}{\kappa}-\frac{\ell}{2}} \sum_{k}  \int_{\mathbb{R}^*}  \frac{ 	|\lambda|^n \,d \lambda}	{ ( |\lambda|(2|k|+n+1)+ 1)^{\frac{m}{2}+\frac{\ell}{4}}} \int_{\mathbb{H}^n }    \frac{ dg}	{ (  |g|^\kappa+1)^{\frac{m}{2}+\frac{\ell}{4}}} \\\nonumber
&\leq C u^{-m+n+1+\frac{2n}{\kappa}+\frac{2}{\kappa}-\frac{\ell}{2}}\sum_{k}  \frac{1}{(2|k|+n+1)^{n+1}}\\&\approx u^{-m+n+1+\frac{2n}{\kappa}+\frac{2}{\kappa}-\frac{\ell}{2}}.
\end{align}
Similarly, we have  
\begin{align}\label{approximation1}
	\left| 	 {tr}\left( Op\left(   (|\lambda|(H+I)+V(g)+u)^{-m}\right) \right)\right| \approx u^{-m+n+1+\frac{2n}{\kappa}+\frac{2}{\kappa}}.
\end{align}
Thus applying trace in (\ref{CH812}) and using (\ref{approximation}), (\ref{approximation1}), we get 
 \begin{align}\label{approximation3}\nonumber
 	\left|\frac{tr(Op(E(g, \lambda, u)))}{tr \left( Op\left( (|\lambda|(H+I)+V(g)+u)^{-m}\right)\right)}\right|&=\left| \frac{ tr\left( (\mathcal{H}+u)^{-m}\right)  }{	 	 {tr}\left( Op\left(   (|\lambda|(H+I)+V(g)+u)^{-m}\right) \right) }-1\right|\\&\leq C\psi_1(u)\to 0~\text{as}~ u\to \infty,
\end{align}
where $\psi_1(u)= \frac{1}{u^{m-1}}+\sum_{\ell=1}^{N} u^{-\frac{\ell}{2}}.$ Note that when $\ell=0, M=0$ and $R_{\ell, M}=1$, (\ref{approximation}), (\ref{approximation1}) has same decay. If $\ell\geq 1$ then (\ref{approximation3}) also holds.
\end{proof}

 Let $w$ be the   complex number varying over the   curve $\Gamma$ (defined in   Section \ref{sec4}). For      fixed $(g, \lambda)\in \mathbb{H}^n\times \mathbb{R}^*$, the class $S _w^m(\mathbb{R}^n)$ defined as
\begin{align}\label{CH80009}\nonumber
&	S_w^m(\mathbb{R}^n)\\&  =\left\{ a\in C^\infty(\mathbb{R}^{2n}\times \Gamma): 	|\partial_{\xi}^{\alpha}\partial_{x}^{\beta} a( x, \xi) | \leq C_{\alpha, \beta}(1+|\lambda|(1+|\xi|^2+|x|^2)+V(g)+|w|)^{\frac{m-|\alpha|}{2}} \right\}.
\end{align}
We obtain the  following result as in  Proposition \ref{CH86}. 
\begin{proposition} \label{New proposition}
	Let $m>0$ be a sufficiently large such that  $	\big({|\lambda|(H+I)+V(g)}+u\big)^{-m}$ is in trace class. Then for  a fixed $(g, \lambda)\in \mathbb{H}^n\times \mathbb{R}^*$, we have 	
	$$	\big({|\lambda|(H+I)+V(g)}+u\big)^{-m}=Op^{W}\left( {\big(|\lambda |(1+|\xi| ^2+|x|^2)+V(g)+u\big)^{-m}}\right)+Op^W(E_{g, \lambda}( u)),$$
	where \begin{align*}
	&	\left| 	tr \big({|\lambda|(H+I)+V(g)}+u\big)^{-m}-tr \left( Op^{W}\left( {\big(|\lambda |(1+|\xi| ^2+|x|^2)+V(g)+u\big)^{-m}}\right)\right) \right|\\&=|tr(Op^W(E_{g, \lambda}( u)))| \\&\leq \psi_2(u) \left|tr \left( Op^{W}\left( {\big(|\lambda |(1+|\xi| ^2+|x|^2)+V(g)+u\big)^{-m}}\right)\right) \right|
	\end{align*}
with  $\psi_2(u)\to 0$ as $u\to \infty.$	 
\end{proposition}
\begin{proof}
	 The proof is based on the same idea as in Proposition \ref{CH86}. For  fixed $(g, \lambda)\in \mathbb{H}^n\times \mathbb{R}^*$, there exists $m\in \mathbb{N}$ such that $	\big({|\lambda|(H+I)+V(g)}+u\big)^{-m}$ is a   trace class operator on $L^2(\mathbb{R}^n)$. We refer to \cite{zel83} for similar pseudo-differential approximation to $\left( |\lambda|(H+I)+V(g)+u\right)^{-m}$ on $L^2(\mathbb{R}^n)$. 	 	 However, we will only indicate some intermediate steps. Now 
	\begin{align}\label{CH8122}
	\big({|\lambda|(H+I)+V(g)}+u\big)^{-m}=Op^{W}\left( {\big(|\lambda |(1+|\xi| ^2+|x|^2)+V(g)+u\big)^{-m}}\right)+Op^W(E_{g, \lambda}( u)),
	\end{align}
	where
	\begin{align*}&E_{g, \lambda}( u)\\&=\sum_{\ell =1}^{N} \sum_{[\frac{\ell}{2}] \leq M\leq \ell }  \frac{\Gamma(s+1)}{\Gamma(s-M)\Gamma(M+1)}   R^{(g, \lambda)}_{\ell, M}(\xi, x ) { ((|\lambda |(1+|\xi| ^2+|x|^2)+V(g)+u\big)^{-m-M}} \\&\qquad +		\frac{1}{2\pi i}\int_{\Gamma} (w+u)^{-m}  S_N^{(g, \lambda)}(w)(|\lambda|(H+I)+V(g)-w)^{-1}  \;d w
	\end{align*}
	with $S_N^{(g, \lambda)}\in S_w^{-N}(\mathbb{R}^n)$.
	  Let $0<s<N$ and $ (I+|\lambda|(H+I)+V(g))^{-\frac{s}{2}}$ is a trace class operator on  $L^2(\mathbb{R}^n)$.  Then imitating the similar  calculations    in \cite{zel83}, we have 
\begin{align*}
	\left| {tr}\left( 	 \int_{\Gamma} (w+u)^{-m} Op^W(S_N^{(g, \lambda)}(w))(|\lambda|(H+I)+V(g)-w)^{-1}  \;d w\right) \right|  
	 \to 0~\text{as}~ u\to \infty,
\end{align*}
and 
\begin{align*}
	&\left|	 {tr}\left( Op^W\left(  R^{(g, \lambda)}_{\ell, M}(\xi, x )  \left(|\lambda |(1+|\xi| ^2+|x|^2)+V(g)+u\right)^{-m-M} \right)\right)\right|\\
	&=\left|\int_{\mathbb{R}^n } \int_{\mathbb{R}^n }   R^{(g, \lambda)}_{\ell, M}(\xi, x )  \left(|\lambda |(1+|\xi| ^2+|x|^2)+V(g)+u\right)^{-m-M}\;dx\;d\xi \right|\\
	&\leq u^{-m-M}\int_{\mathbb{R}^n } \int_{\mathbb{R}^n }  \left| R^{(g, \lambda)}_{\ell, M}(\xi, x ) \right|  \left(u^{-1}|\lambda |(1+|\xi| ^2+|x|^2)+u^{-1}V(g)+1\right)^{-m-M}\;dx\;d\xi\\
		&\leq u^{-m-M}\int_{\mathbb{R}^n } \int_{\mathbb{R}^n }  \left| R^{(g, \lambda)}_{\ell, M}(\xi, x ) \right|  \left(u^{-1}|\lambda |(|\xi| ^2+|x|^2)+1\right)^{-m-M}\;dx\;d\xi\\
	&=u^{-m-M+ {n}}\int_{\mathbb{R}^n } \int_{\mathbb{R}^n }  \left| R^{(g, \lambda)}_{\ell, M}(u^{\frac{1}{2}}\xi, u^{\frac{1}{2}}x ) \right|  \left(|\lambda |(|\xi| ^2+|x|^2)+1\right)^{-m-M}\;dx\;d\xi\\
	&\leq u^{-m-M+ {n}}\int_{\mathbb{R}^n } \int_{\mathbb{R}^n }  \left| (1+|\lambda|(1+u|\xi|^2+u|x|^2)+V(g)+|w|)^{\frac{2M-\ell}{2}}	 \right| \\&\qquad \times  \left(|\lambda |(|\xi| ^2+|x|^2)+1\right)^{-m-M}\;dx\;d\xi\\
		&\leq C u^{-m+n-\frac{\ell}{2}}.
\end{align*}
Note that if $\ell=1, M=1$ and $ R^{(g, \lambda)}_{1, 1}(x, \xi)=1$. So for $\ell\geq 2,$ 
$$\left|	 {tr}\left( Op^W\left(  R^{(g, \lambda)}_{\ell, M}(\xi, x )  \left(|\lambda |(1+|\xi| ^2+|x|^2)+V(g)+u\right)^{-m-M} \right)\right)\right|\approx u^{-m+n-1}.$$
Similarly, we have  
\begin{align*}
	&\left|	 {tr}\left( Op^W\left(    \left(|\lambda |(1+|\xi| ^2+|x|^2)+V(g)+u\right)^{-m-M} \right)\right)\right|\approx   u^{-m+n }.
\end{align*}
Thus 
\begin{align*}
	&\left| \frac{ tr(Op^W(E_{g, \lambda}( u)))}{ tr \left( Op^{W}\left( {\big(|\lambda |(1+|\xi| ^2+|x|^2)+V(g)+u\big)^{-m}}\right)\right) }\right|\\&=\left| \frac{	 	 {tr}\left( Op\left(   (|\lambda|(H+I)+V(g)+u)^{-m}\right) \right) }{ tr\left( Op^{W}\left( {\big(|\lambda |(1+|\xi| ^2+|x|^2)+V(g)+u\big)^{-m}}\right)\right)  }-1\right|\leq \psi_2(u) \to 0~\text{as}~ u\to \infty,
\end{align*}
where $\psi_2(u)= \frac{1}{u^{m-1}}+\sum_{\ell=1}^{N} u^{-1}.$
\end{proof}
\begin{remark}\label{DCT}
Note that  for sufficiently large $m\in \mathbb{N}$, the operator $$Op  \left( Op^W\left(  R^{(g, \lambda)}_{\ell, M}(\xi, x )  \left(|\lambda |(1+|\xi| ^2+|x|^2)+V(g)+u\right)^{-m-M} \right)\right)$$ is a trace class operator  on $L^2(\mathbb{H}^n),$ since from    Proposition \ref{New proposition}, we have  
\begin{align*} 
	& \left|tr\left(Op  \left( Op^W\left(  R^{(g, \lambda)}_{\ell, M}(\xi, x )  \left(|\lambda |(1+|\xi| ^2+|x|^2)+V(g)+u\right)^{-m-M} \right)\right)\right)\right|\\ 
	&\leq \int_{\mathbb{H}^n } 
	\int_{\mathbb{R}^*}\left|	 {tr}\left( Op^W\left(  R^{(g, \lambda)}_{\ell, M}(\xi, x )  \left(|\lambda |(1+|\xi| ^2+|x|^2)+V(g)+1\right)^{ {-m-M}} \right)\right)\right|dg\;d\mu(\lambda)\\ 
&\leq  \int_{\mathbb{H}^n } 
\int_{\mathbb{R}^*} \int_{\mathbb{R}^{2n} }    \left| (|\lambda|(1+|\xi|^2+|x|^2)+V(g)+1)^{\frac{2M-\ell}{2}}	 \right|  \\&\qquad\times   \left(|\lambda |(1+|\xi| ^2+|x|^2)+V(g)+1\right)^{-m-M}\;dx\;d\xi\;dg\;d\mu(\lambda)<\infty.\\
&\leq  \int_{\mathbb{H}^n } 
\int_{\mathbb{R}^*} \int_{\mathbb{R}^{2n} }     (1+|\lambda|(1+|\xi|^2+|x|^2)+V(g))^{\frac{-m-\ell}{2}}	    \;dx\;d\xi\;dg\;d\mu(\lambda)\\
&\leq C\int_{\mathbb{R}^{2n} }     \int_{\mathbb{H}^n } \int_{\mathbb{R}^*}  \frac{\lambda^n}{(1+|\lambda|(1+|\xi|^2+|x|^2)+V(g)+1)^{\frac{m+\ell}{2}}}\;dx\;d\xi \,d g \,d\lambda \\
&=C\int_{\mathbb{R}^{2n} } \frac{1}{(1+|\xi|^2+|x|^2)^{{n+1}}  } \;dx\;d\xi \int_{\mathbb{H}^n } 	\int_{1+V(g)}^\infty    \frac{(u-1-V(g))^n}{u^{\frac{m+\ell}{2}}} \, du\,d g   \\
&=C\int_{\mathbb{R}^{2n} }  \frac{1}{(1+|\xi|^2+|x|^2)^{ {n+1}}  }  \;dx\;d\xi \int_{\mathbb{H}^n } \frac{1}{(1+V(g))^{\frac{m+\ell}{2}-n-1}} \,d g   <\infty.
\end{align*}
Similarly it can be shown that $ Op  \left( Op^W\left(   \left(|\lambda |(1+|\xi| ^2+|x|^2)+V(g)+u\right)^{-m} \right)\right) $ is a trace class operator on $L^2(\mathbb{H}^n)$ for sufficiently large $m\in \mathbb{N}$.
\end{remark}
We take the poitive integer $m$ such that the requirement
of $m$-th power of the operators discussed earlier  to be  a trace class operator. 
 \section{Szeg\"o limit theorem for $\mathcal{H}$}\label{sec5}
Now we are in a position to prove our main results. We start this section with  the following lemmas.
\begin{lemma}
	Let $M_{\textbf{b}}$ be the multiplication operator defined in Theorem \ref{mo}, then  for any $f\in C(\mathbb{R})$, ${tr}{f(\mathcal{P}_r M_\textbf{b} \mathcal{P}_r)} = {tr}{(\mathcal{P}_rf(M_{\textbf{b}})\mathcal{P}_r)}$
\end{lemma}
\begin{proof}
Notice that  $\| (I-\mathcal{P}_r)M_{\textbf{b}} \mathcal{P}_r\|_{B_2}^2= {tr}{(\mathcal{P}_r M_\textbf{b} \mathcal{P}_r)}= {tr}{(\mathcal{P}_rM_{\textbf{b}}^2 \mathcal{P}_r)}-{tr}{(\mathcal{P}_rM_{\textbf{b}}\mathcal{P}_r)^2}.$ Also $\mathcal{P}_rM_{\textbf{b}}^2 \mathcal{P}_r$ is an operators on $L^2(\mathbb{H}^n)$ with kernel $K_1(g,g_1)=\displaystyle\sum_{k_1,k_2\leq r}\langle \textbf{b}^2 e_{k_1},e_{k_2}\rangle e_{k_2}(g)\overline{e_{k_1}(g_1)}$, for any orthonormal basis $\{e_k\}$ of $L^2(\mathbb{H}^n)$. Therefore ${tr}{(\mathcal{P}_rM_{\textbf{b}}^2 \mathcal{P}_r)} =\int_{\mathbb{H}^n}K_1(g,g)dg=\displaystyle\sum_{k\leq r}\langle \textbf{b}^2 e_{k},e_{k}\rangle$.
Further, ${tr}{(\mathcal{P}_rM_{\textbf{b}}\mathcal{P}_r)^2}=\int_{\mathbb{H}^n}K_2(g,g)dg=\displaystyle\sum_{k\leq r}\langle \textbf{b}^2 e_{k},e_{k}\rangle,$ where the operator $\mathcal{P}_rM_{\textbf{b}}\mathcal{P}_rM_{\textbf{b}}\mathcal{P}_r$  is an integral operator with kernel $$K_2(g,g_1)=\displaystyle\sum_{k_1,k_2,k_3\leq r}\langle \textbf{b} e_{k_1},e_{k_2}\rangle \langle \textbf{b} e_{k_2},e_{k_3}\rangle e_{k_3}(g)\overline{e_{k_1}(g_1)}.$$ So ${tr}{(\mathcal{P}_rM_{\textbf{b}}^2 \mathcal{P}_r)} = {tr}{(\mathcal{P}_rM_{\textbf{b}}\mathcal{P}_r)^2}.$ So $\| (I-\mathcal{P}_r)M_{\textbf{b}} \mathcal{P}_r\|_{B_2}^2=0$.  Observe that  for $n\in \mathbb{N}$,  $\mathcal{P}_rM_{\textbf{b}}^n \mathcal{P}_r = \mathcal{P}_rM_{\textbf{b}} (\mathcal{P}_r+(I-\mathcal{P}_r))M_{\textbf{b}} \cdots M_{\textbf{b}} \mathcal{P}_r = (\mathcal{P}_rM_{\textbf{b}} \mathcal{P}_r)^n$+ terms with a factor of $(I-\mathcal{P}_r)M_{\textbf{b}} \mathcal{P}_r$. By Cauchy-Schwarz inequality, $tr (\mbox{terms with a factor of}$ $(I-\mathcal{P}_r)M_{\textbf{b}})$  is dominated by some constant (depending on $\textbf{b}$) times
  $\| (I-\mathcal{P}_r)M_{\textbf{b}} \mathcal{P}_r\|_{B_2}$. Therefore $|{tr}{(\mathcal{P}_rM_{\textbf{b}}^n \mathcal{P}_r)}-{tr}{(\mathcal{P}_rM_{\textbf{b}}\mathcal{P}_r)^n}|=0$.
  Thus ${tr}{f(\mathcal{P}_rM_{\textbf{b}}\mathcal{P}_r)} = {tr}{(\mathcal{P}_rf(M_{\textbf{b}})\mathcal{P}_r)}$ for $f(x)=x^n$, $\forall n \in \mathbb{N}$ and this result can be extended to continuous functions as an application of the Weierstrass approximation theorem and spectral theorem.
  \end{proof}
\begin{lemma}\label{piir}
For $r>0$ define $I_{r}: L^2(\mathbb{H}^n) \rightarrow L^2(\mathbb{H}^n)$ by $$I_{r}(\phi)(g) = \displaystyle\int_{-r}^{r}tr(\pi_\lambda^*(g)\hat{\phi}(\lambda))d\mu(\lambda).$$ Then
\begin{align}\label{aa}
\lim_{r\to\infty} \frac{{tr}{(\mathcal{P}_rf(M_{\textbf{b}})\mathcal{P}_r)}}{tr~(\mathcal{P}_r)} &= \lim_{r\to\infty} \frac{{tr}{(\mathcal{P}_rI_{r}f(M_{\textbf{b}})\mathcal{P}_rI_{r})}}{tr~(\mathcal{P}_rI_{r})}.
\end{align}
\end{lemma}
\begin{proof}
	We know that if $X$ is a positive trace class operator and $Y$ is a bounded operator on $L^2(\mathbb{H}^n)$ then $|tr~(XYX)|\leq \|Y\| |tr~(X^2).$ Using this inequality  we get
	$ \big|tr~(\mathcal{P}_r) - tr~(\mathcal{P}_rI_{r})\big |= \big |tr~(\mathcal{P}_r(I-I_{r}))\big |\leq \|I-I_{r}\| |tr~(\mathcal{P}_r)|.$
But for $\psi \in L^2(\mathbb{H}^n) $, an application of Plancherel formula gives $\|(I-I_{r})\psi \|_2^2=\displaystyle\int_{|\lambda| >r}\|\hat{\psi }(\lambda)\|_{B_2}^2 \,d\mu(\lambda)\rightarrow 0  ~ \mathrm{as} ~ r \rightarrow \infty.$
		Therefore, 	\begin{eqnarray}\label{a}	\bigg|\frac{tr~(\mathcal{P}_rI_{r})}{tr~(\mathcal{P}_r)} - 1\bigg| \leq \|I-I_{r}\| \rightarrow 0  ~ \mathrm{as} ~ r \rightarrow \infty.\end{eqnarray}

We add a suitable constant to make the operator $M_{\textbf{b}}$ positive and any $f\in C(\mathbb{R})$ can be written as the difference of two positive functions namely the positive and the negative part of $f$. So without loss of generality we take $f(M_{\textbf{b}})$ as a positive operator.
Further,
\begin{align*}
&\big |{tr}{(\mathcal{P}_rf(M_{\textbf{b}})\mathcal{P}_r)- {tr}(\mathcal{P}_rI_{r}f(M_{\textbf{b}})\mathcal{P}_rI_{r})}\big | \\
&= \big |{tr}(\mathcal{P}_rf(M_{\textbf{b}})\mathcal{P}_r(I-I_{r}))- {tr}{(\mathcal{P}_r(I-I_{r})f(M_{\textbf{b}})\mathcal{P}_rI_{r})}\big | \\
& \leq \big |{tr}(\mathcal{P}_rf(M_{\textbf{b}})\mathcal{P}_r)\big | \|I-I_{r}\| + \big |{tr}{(\mathcal{P}_rf(M_{\textbf{b}})\mathcal{P}_rI_{r})}\big | \|I-I_{r}\|.
\end{align*}
Therefore,
\begin{eqnarray}\label{b}
\bigg|\frac{{tr}(\mathcal{P}_rI_{r}f(M_{\textbf{b}})\mathcal{P}_rI_{r})}{{tr}(\mathcal{P}_rf(M_{\textbf{b}})\mathcal{P}_r)}-1\bigg|
 \leq  \left(1+ \bigg  |\frac{{tr}{(\mathcal{P}_rf(M_{\textbf{b}})\mathcal{P}_rI_{r})}}{{tr}(\mathcal{P}_rf(M_{\textbf{b}})\mathcal{P}_r)}\bigg|\right)\|I-I_{r}\|
  \rightarrow 0,
\end{eqnarray} as $r \rightarrow \infty$.
Combining (\ref{a}) and (\ref{b}) we get (\ref{aa}).
\end{proof}

\noindent{\bf Proof of theorem \ref{mo}:}\\
The operator $\mathcal{P}_rI_{r}f(M_{\textbf{b}})\mathcal{P}_rI_{r}$ is an integral operator with kernel $$K_{r}(g, g_1)=\int_{-r}^{r} {tr}(\pi_{\lambda}^*(g)I_{r\times r}f(\textbf{b} (g_1))\pi_{\lambda}(g_1) )\,d \mu(\lambda).$$
Therefore
\begin{align*}
{tr}{(\mathcal{P}_rI_{r}f(M_{\textbf{b}})\mathcal{P}_rI_{r})}
=\int_{\mathbb{H}^n} K_{r}(g, g)\,dg
&= r\int_{-r}^{r}  d \mu(\lambda) \int_{\mathbb{H}^n}f( \textbf{b} (g)) \,dg
\end{align*}
and
$$
{tr}{(\mathcal{P}_rI_{r})}=\displaystyle\sum_{i\leq r}\langle \widehat{I_{r} \phi_i}, \hat{\phi_i}\rangle
=\sum_{i\leq r}\langle \hat{ \phi_i}\mathcal{X}_{[-r, r]}, \hat{\phi_i}\rangle
=\sum_{i\leq r}\int_{-r}^{r} {tr} (\hat{ \phi_i}^*(\lambda) \hat{\phi_i}(\lambda)) \,d\mu(\lambda)
= r\int_{-r}^{r} \,d\mu(\lambda),
$$
for any orthonormal basis $\{\phi_i\}_{i=1}^\infty$ of  $L^2(\mathbb{H}^n).$
Thus \begin{align*}
\lim_{r\to\infty} \frac{{tr}{(\mathcal{P}_rf(M_{\textbf{b}})\mathcal{P}_r)}}{tr~(\mathcal{P}_r)} &= \lim_{r\to\infty} \frac{{tr}{(\mathcal{P}_rI_{r}f(M_{\textbf{b}})\mathcal{P}_rI_{r})}}{tr~(\mathcal{P}_rI_{r})}
= \int_{\mathbb{H}^n}f(\textbf{b}(g)) \,dg.
\end{align*}\vspace{0.5cm}

\noindent{\bf Proof of theorem \ref{sch}:}\\
 Now we prove Szeg\"o limit theorem for $\mathcal{H}$ under certain assumptions on the symbol $a(g,\lambda)$ to ensure the existence the RHS limit in Theorem \ref{sch} (see \cite{zel83}).  We assume \begin{eqnarray}\label{assum}\displaystyle\lim_{E\to\infty} \bar{a}(E)=a,\end{eqnarray} where $\bar{a}(E)=\dfrac{1}{S(E)}\displaystyle \int_{G_E}  a_{g, \lambda }(\xi, x) \,d\xi\,dx \, dg \, d\mu(\lambda) $ and  $S(E)=\displaystyle\displaystyle \int_{G_E} \,d\xi\,dx \, dg \, d\mu(\lambda), $ with $G_{E}=\{(g, \lambda, \xi, x)\in \mathbb{H}^n \times \mathbb{R}^*\times  \mathbb{R}^n\times \mathbb{R}^n : |\lambda |(1+|\xi| ^2+|x|^2)+V(g)= E \}$ and \begin{eqnarray}\label{assum1}V(g)\sim V_0|g|^\kappa ~~\mbox{as}~|g|\to \infty \end{eqnarray}~~for~ real~~ $\kappa>0$ in~the~sense~that $V(g)=V_0|g|^\kappa+W(g),~\mbox{where}~W(g)=o(|g|^\kappa)$.
\begin{proposition}\label{assymptotic}
Let $G^E=\{(g, \lambda, \xi, x)\in \mathbb{H}^n \times \mathbb{R}^*\times  \mathbb{R}^n\times \mathbb{R}^n : |\lambda |(1+|\xi| ^2+|x|^2)+V(g)\leq E \}$. Then volume of $G^{E}=v(E) \approx   E^{n+1+\frac{2(n+1)}{\kappa}}$  as $E \rightarrow \infty $.
\end{proposition}
\begin{proof} : Using  the homogeneous norm on $\mathbb{H}^n$ we have
\begin{align*}
v(E) &=\displaystyle \int_{G^{E}}  \,d\xi\,dx \, dg \, d\mu(\lambda) \\
&=C_n\displaystyle  \int_{\mathbb{H}^n} \iint_{\mathbb{R}^{2n}} \left ( \int_{|\lambda| \leq \frac{(E-V(g))_{+}}{ 1+|x|^2+|\xi|^2}}   |\lambda|^n \, d\lambda \right) \,d\xi\,dx   \, dg \\
&=2C_n\displaystyle  \int_{\mathbb{H}^n}(E-V(g))_{+}^{n+1}  \, dg  \iint_{\mathbb{R}^{2n}}    \left( \frac{1}{ 1+|x|^2+|\xi|^2} \right)^{n+1}  \,d\xi\,dx  \\
&=C_n' E^{n+1}\displaystyle  \int_{V\leq E} \left ( {1-E^{-1}V(g)}\right)^{n+1} \, dg\\
&\sim  C_n'E^{n+1}\displaystyle  \int_{V\leq E} \left ( {1-E^{-1}(V_0|g|^\kappa+W(g))}\right)^{n+1} \, dg\\
&=C_n' E^{n+1+\frac{2(n+1)}{\kappa}}\displaystyle  \int_{V\leq E} \left ( 1-V_0|g|^k-E^{-1}W(\widetilde{g})\right)^{n+1} \, dg,
\end{align*} where $\widetilde{g}=(E^\frac{1}{\kappa}x_1, E^\frac{1}{\kappa}x_2,\cdots, E^\frac{1}{\kappa}x_{2n},E^\frac{2}{\kappa}t)$ for $g=(x_1,x_2\cdots, x_{2n},t)\in\mathbb{H}^n.$ Since $ \lim_{E\to\infty}$ $E^{-1}W(\widetilde{g})=0$, the right hand side of the above integral converges to $\int_{V\leq E} \left ( 1-V_0|g|^\kappa\right)^{n+1} \, dg$ by dominated convergence theorem.
\end{proof}
\begin{lemma}\label{ratio}
Let  $\phi(r)={tr}({\mathcal{P}_r})$ and $\psi(r)=tr({\mathcal{P}_rA\mathcal{P}_r})$. Then under the assumption (\ref{assum}) and (\ref{assum1}), we have
$$\Phi(u)=\int_{0}^{\infty} \frac{\phi(r)}{(r+u)^{m+1}} \, dr=  \int_{0}^{\infty}\frac{\displaystyle \int_{G^E} \,d\xi\,dx \, dg \, d\mu(\lambda)}{(E+u)^{m+1}}\,dE+E_1(u)\\$$ and $$\Psi(u)=\int_{0}^{\infty} \frac{\psi(r)}{(r+u)^{m+1}} \,dr= \int_{0}^{\infty}\frac{\displaystyle \int_{G^E} a_{g, \lambda}(\xi, x)
	\,d\xi\,dx \, dg \, d\mu(\lambda)}{(E+u)^{m+1}}\,dE+E_2(u),$$
 with $|E_i(u)| \to0$ as $u\to\infty,~i=1,2.$
\end{lemma}
\begin{proof}
	The operator $\mathcal{H}$ has discrete spectrum of eigenvalues $0\leq c_1\leq c_2\cdots \infty$. Let $\{\psi_j\}_{j=1}^\infty$ be the complete set of eigenfunctions on corresponding to the eigenvalues $\{c_j\}$  on  $L^2(\mathbb{H}^n)$.
 Then $
\psi(r)= tr ({\mathcal{P}_r}A{\mathcal{P}_r})= \displaystyle\sum_{c_j\leq r} \langle A\psi_j  ,\psi_j \rangle
$  and
$\psi'(r)= \displaystyle\sum_{j=1}^{\infty} \langle A\psi_j  ,\psi_j \rangle \delta(r-c_j).$
Now 
\begin{align*}
\Psi(u)&= \int_{0}^{\infty} \frac{\psi(r)}{(r+u)^{m+1}} \,dr= m\sum_{j=1}^{\infty} \langle A \psi_j  ,\psi_j \rangle~\int_{0}^{\infty}\frac{\delta(r-c_j)}{(r+u)^{m}}\, dr\\&= m\sum_{j=1}^{\infty} \langle A\psi_j  ,\psi_j \rangle~\frac{1}{(c_j+u)^{m}}=m~{tr}\left({A{(\mathcal{H}+u)^{-m}}}\right).
\end{align*}
By Proposition \ref{CH86}, we obtain 
\begin{eqnarray*}\Psi(u)&=&m~{tr}\left({A{(\mathcal{H}+u)^{-m}}}\right)\\&=&m~{tr}\left({A~Op\big({|\lambda|(H+I)+V(g)}+u\big)^{-m}}\right) +m~{tr}\left(A~Op(E(g, \lambda, u)) \right)
\end{eqnarray*} 
with  $|{tr}\left(A~Op(E(g, \lambda, u)) \right)|\leq \|A\|\left|  {tr}\left(Op(E(g, \lambda, u)) \right)\right| \to 0$ $u\to \infty$. Thus  for large $u$, using  (\ref{schmidt}) and (\ref{hilbert1}), we have 
\begin{align*}&\Psi(u)=m\int_{\mathbb{H}^n\times \mathbb{R}^*} {tr}\left({{a(g, \lambda)\big({|\lambda|(H+I)+V(g)}+u\big)^{-m}}}\right)\, dg \, d\mu(\lambda)\\
	&=m\int_{\mathbb{H}^n\times \mathbb{R}^*} {tr} 
	\bigg(a(g, \lambda) Op^{W}\left( {\big(|\lambda |(1+|\xi| ^2+|x|^2)+V(g)+u\big)^{-m}} \right)\bigg.\\&\bigg. \qquad +a(g, \lambda)Op^W(E_{g,\lambda}(u))\bigg) \, dg \, d\mu(\lambda)\\
	&=m \int_{\mathbb{H}^n\times \mathbb{R}^*}\int_{\mathbb{R}^n\times \mathbb{R}^n} \left( a_{g\lambda}(\xi, x){\big(|\lambda |(1+|\xi| ^2+|x|^2)+V(g)+u\big)^{-m}}\right)\,d\xi\,dx\, dg \, d\mu(\lambda) +E_1(u)\\
	&= \int_{0}^{\infty}\frac{\displaystyle \int_{G_E} a_{g\lambda}(\xi, x)\,d\xi\,dx \, dg \, d\mu(\lambda)}{(E+u)^{m+1}}\,dE+E_1(u),
\end{align*} 
where $E_1(u)=m \displaystyle\int_{\mathbb{H}^n\times \mathbb{R}^*}{tr}(a(g, \lambda)Op^W(E_{g,\lambda}(u)))~dg d\mu(\lambda)$. From Remark \ref{DCT}, we   conclude that $|E_1(u)|\to 0$ as $u\to \infty$ by   dominated convergence theorem.
Similarly taking $A=I$, we get  $\phi(r)={tr}{(\mathcal{P}_r})$, and in this case, for large $u$, we have
\begin{align*}
\Phi(u)&= \int_{0}^{\infty}\frac{\displaystyle \int_{G_E}  
		 \,d\xi\,dx \, dg \, d\mu(\lambda)}{(E+u)^{m+1}}\,dE +E_2(u)
\end{align*}
with $|E_2(u)|\to 0$ as $u\to\infty.$
  %
%
%
%
%
\end{proof}
In order to prove the Szeg\"o limit theorem for the Schr\"odinger operator $\mathcal{H}$,  we need to estimate the asymptotic growth of the measures $tr(\mathcal{P}_rA\mathcal{P}_r)$ and $tr(\mathcal{P}_r)$.  Therefore we apply Keldysh Tauberian Theorem (see Theorem 5.4 in Appendix) to compare the measures.
\begin{corollary}\label{abs}
	Consider the  self-adjoint operator $\mathcal{P}_r$ and $v(r)$  as given in Theorem \ref{sch} and Proposition \ref{assymptotic} respectively. Let $\phi(r)=tr (\mathcal{P}_r), $ $ \psi(r)=tr (\mathcal{P}_rA\mathcal{P}_r)$ then we have the following asymptotic :
	\begin{enumerate}
		\item $v (r) \approx C r^{n+1+\frac{2(n+1)}{\kappa}}$ as  $r \to \infty.$
		\item $v  $ is multiplicatively continuous.
		\item $tr (\mathcal{P}_r) \approx C r^{n+1+\frac{2(n+1)}{\kappa}}$ as  $r \to \infty.$
		\item $\displaystyle\sup_{\mu\leq r}[tr(\mathcal{P}_{\mu+r_1})-tr(\mathcal{P}_\mu)] \leq  tr (\mathcal{P}_r)\left[   \left({n+1+\frac{2(n+1)}{\kappa}}\right)  \frac{r_1}{r}+  \mathcal{O}\left(\frac{1}{r}\right)^2 \right]$, as $r \to \infty.$
		\item $\psi $  is  multiplicatively continuous.
	\end{enumerate}
\end{corollary}
\begin{proof}
	Clearly (1)  directly follows from  Proposition \ref{assymptotic}. Now
	\begin{eqnarray*}
		\displaystyle \lim_{\substack {r \to \infty\\  \tau \to 1}} \frac{v(\tau r)}{v(r)}= \displaystyle \lim_{\substack {r \to \infty\\  \tau \to 1}} \frac{(\tau r)^{{n+1+\frac{2(n+1)}{\kappa}}}} {r^{{n+1}+\frac{2(n+1)}{\kappa}}}=\displaystyle \lim_{ \tau \to 1}\tau {^{{n+1}+\frac{2(n+1)}{\kappa}}}=1.
	\end{eqnarray*}
	Therefore $v $ is multiplicatively continuous function. We choose sufficiently large $m$ such that the operator $(\mathcal{H}+uI)^{-m}$ is a trace class operator. Therefore by Lemma \ref{ratio} and Theorem 8  of Grishin-Poedintseva \cite{gri}, we get
	$	\phi (r)/  v(r) \to 1 ~ \mbox{as} ~r \to \infty. $
	This proves (3).
	Using the asymptotic in (3), it is easy to check that
	\begin{align*}
	\sup_{\mu\leq r}[tr(\mathcal{P}_{\mu+r_1})-tr(\mathcal{P}_\mu)] \leq  tr(\mathcal{P}_r)\left[ \left(n+1+\frac{2(n+1)}{\kappa}\right) {{\frac{r_1}{r}+  \mathcal{O}\left(\frac{1}{r}\right)^2}} \right].
	\end{align*}
	To prove (5), notice that if $\varphi$ and $\chi$ are two distribution functions satisfying $\displaystyle \lim_{r \to \infty } \frac{\varphi(r)}{\chi(r)}=1$ then $\varphi$ is multiplicatively continuous whenever $\chi$ is.  Therefore $\psi  $ is   also a multiplicatively continuous function.
\end{proof}
\begin{theorem} Under the assumption (\ref{assum}) and (\ref{assum1}) we have
\begin{align*}
\lim_{r\to\infty} \frac{{tr}{(\mathcal{P}_rA\mathcal{P}_r)}}{tr~(\mathcal{P}_r)} &= \lim_{r\to\infty}  \frac{\int_{G^{r}} a_{g, \lambda}(\xi, x)
	\,d\xi\,dx \, dg \, d\mu(\lambda) }{\int_{G^{r}}
\,d\xi\,dx \, dg \, d\mu(\lambda) }.
\end{align*}
\end{theorem}
\begin{proof}
	The proof follows directly by Lemma \ref{ratio}, as all the requirements (by our assumption (\ref{assum}) on the symbol $a(g,\lambda)$) of Theorem 8  of Grishin-Poedintseva \cite{gri} are satisfied.
\end{proof}
\begin{corollary}\label{final}
	Let $P(r)$ be a polynomial in $\mathbb{R}.$ Then
	\begin{align*}
	\lim_{r\to\infty} \frac{{tr}({\mathcal{P}_rP(A)\mathcal{P}_r})}{tr~(\mathcal{P}_r)} &=  \lim_{r\to\infty}  \frac{\int_{G^{r}} P(a_{g, \lambda}(\xi, x))
		\,d\xi\,dx \, dg \, d\mu(\lambda)}{\int_{G^{r}}
		\,d\xi\,dx \, dg \, d\mu(\lambda) }
	\end{align*}
\end{corollary}
	\begin{proof}
	From  the asymptotic expression of Theorem \ref{hcom} along with Remark \ref{ab},  we   see that 	the operator $P(A)$ has symbol $P(a(g, \lambda ))+E(g, \lambda)+E_{-1}(g, \lambda )$, where $E(g, \lambda), E_{-1}(g, \lambda )\in S_\mathcal{H}^{-1}(\mathbb{H}^n)$ (the term associated with  $[\alpha]=1$ is $E$ and $E_{-1} $  is the remaining terms   with $[\alpha]>1$  in the asymptotic expansion).   The proof will be complete if we   show   
\begin{align}\label{associated}
	 \lim_{r\to\infty}  \frac{\int_{G^{r}} \tilde{E}_{g, \lambda}(\xi, x)
		\,d\xi\,dx \, dg \, d\mu(\lambda)}{\int_{G^{r}}
		\,d\xi\,dx \, dg \, d\mu(\lambda) }=0, 
\end{align}
	where $E(g, \lambda)+E_{-1}(g, \lambda )=Op^W(\tilde{E}_{g, \lambda}(\xi, x))$. Now 
proceeding   as in  proposition \ref{assymptotic}, we get 
	\begin{align*}
&	 \int_{G^{r}}| \tilde{E}_{g, \lambda}(\xi, x)| \,d\xi\,dx \, dg \, d\mu(\lambda) \\
	 		&\leq C  \int_{\mathbb{H}^n} \iint_{\mathbb{R}^{2n}} \left ( \int_{|\lambda| \leq \frac{(r-V(g))_{+}}{ 1+|x|^2+|\xi|^2}}  (1+|\lambda|(1+|\xi|^2+|x|^2)+V(g))^{-\frac{1}{2}} |\lambda|^n \, d\lambda \right) \,d\xi\,dx   \, dg \\
		&\leq  C  \int_{\mathbb{H}^n} \iint_{\mathbb{R}^{2n}} \left ( \int_{|\lambda| \leq \frac{(r-V(g))_{+}}{ 1+|x|^2+|\xi|^2}}   |\lambda|^{n-\frac{1}{2}} \, d\lambda \right) \,d\xi\,dx   \, dg \\
		&=C  \int_{\mathbb{H}^n}(r-V(g))_{+}^{n+\frac{1}{2}}  \, dg  \iint_{\mathbb{R}^{2n}}    \left( \frac{1}{ 1+|x|^2+|\xi|^2} \right)^{n+\frac{1}{2}}  \,d\xi\,dx  \\
		&=C r^{n+\frac{1}{2}}\displaystyle  \int_{V\leq r} \left ( {1-r^{-1}V(g)}\right)^{n+\frac{1}{2}} \, dg\\
		&\approx C r^{n+\frac{1}{2}+\frac{2(n+1)}{\kappa}}.
	\end{align*} 
On the other hand, from Corollary \ref{abs}, we have   $tr (\mathcal{P}_r) \approx r^{n+1+\frac{2(n+1)}{\kappa}}$.  Thus we get (\ref{associated}).
	\end{proof}
\begin{lemma}\label{direct}
	Let $\mathcal{H}, A$  be the operators defined in Theorem \ref{sch}.  Then

(a) $\sqrt{\mathcal{H}}=\mathcal{H}_{\frac{1}{2}}+ C$, where  $\mathcal{H}_{\frac{1}{2}}=Op( {H}_{\frac{1}{2}}(g, \lambda))$ and $C$ is a bounded operator on $L^2(\mathbb{H}^n)$,

	(b)  the operator $[\sqrt{\mathcal{H}}, A]$ is bounded on $ L^2(\mathbb{H}^n)$,
	
	(c) under the 	assumptions  of Theorem \ref{sch}, we have
	$$\bigg | \frac{{tr}{f(\mathcal{P}_rA\mathcal{P}_r)}}{tr~(\mathcal{P}_r)}  - \frac{{tr}{(\mathcal{P}_rf(A)\mathcal{P}_r)}}{tr~(\mathcal{P}_r)} \bigg |\rightarrow 0 ~\text{as} ~ r \rightarrow \infty.$$
\end{lemma}
\begin{proof}

(a) Let $f(w)=w^{\frac{1}{2}}$. Proceeding as in Subsection \ref{subsec4.1}, we get   $\sqrt{\mathcal{H}}=\mathcal{H}_{\frac{1}{2}} + F_{\frac{1}{2}}$, where $\mathcal{H}_{\frac{1}{2}}=Op( {H}_{\frac{1}{2}}(g, \lambda))$ with ${H}_{\frac{1}{2}}(g, \lambda) \in S_\mathcal{H}^{1}(\mathbb{H}^n)$   and $F_{\frac{1}{2}}$ is defined in (\ref{errorterm}) with $S_N \in S_{\mathcal{H},w}^{-N}(\mathbb{H}^n).$ We choose $N>0$ such that the integral  (\ref{errorterm})  converges in the norm on $B(L^2(\mathbb{R}^n))$. Denoting  $C=F_{\frac{1}{2}}$, we have $\sqrt{\mathcal{H}}=\mathcal{H}_{\frac{1}{2}} + C$ as desired. 

(b) From part (a) we have $\sqrt{\mathcal{H}}=\mathcal{H}_{\frac{1}{2}} + C$.
 Since $C$ is bounded,  $[\sqrt{\mathcal{H}}, A]$ is bounded if  $[\mathcal{H}_{\frac{1}{2}}, A]$ is bounded on $ L^2(\mathbb{H}^n)$. Now using the composition formula  (\ref{com}) of Theorem \ref{hcom}, there exist  two symbols $R_1(g, \lambda), R_2(g, \lambda)\in S_\mathcal{H}^0(\mathbb{H}^n) $ such that 
\begin{align}\label{square}\nonumber
&	[\mathcal{H}_{\frac{1}{2}}, A]=\mathcal{H}_{\frac{1}{2}} A- A \mathcal{H}_{\frac{1}{2}} \\\nonumber
&= Op\left( {H}_{\frac{1}{2}}(g, \lambda)\#_{\mathbb{H}^n} a(g, \lambda )\right)-Op\left(a(g, \lambda ) \#_{\mathbb{H}^n} {H}_{\frac{1}{2}}(g, \lambda)\right)\\\nonumber
&= Op\left(\ {H}_{\frac{1}{2}}(g, \lambda)a(g, \lambda)+\Delta  {H}_{\frac{1}{2}}(g, \lambda)~X_ga(g, \lambda )+   R_1(g, \lambda)\right )\\\nonumber
&\qquad -  Op\left(a(g, \lambda )  {H}_{\frac{1}{2}}(g, \lambda) +\Delta a(g, \lambda ) X_g {H}_{\frac{1}{2}}(g, \lambda)+   R_2(g, \lambda) \right)\\\nonumber
&= Op\left(Op^W(a_{g, \lambda}(\xi,u)  \#    {H}_{\frac{1}{2},g, \lambda}(\xi, u)-  {H}_{\frac{1}{2},g, \lambda}(\xi, u)\# a_{g, \lambda}(\xi,u)\right  ) \\\nonumber
&\qquad+Op\left( \Delta{H}_{\frac{1}{2}}(g, \lambda)~X_ga(g, \lambda )-\Delta a(g, \lambda ) X_g {H}_{\frac{1}{2}}(g, \lambda)\right )+  Op\left(   R_1(g, \lambda)-R_2(g, \lambda)\right )\\\nonumber
&= Op\left(Op^W(F_{g, \lambda}^{1}(\xi,u)  +F_{g, \lambda}^{2}(\xi,u))\right) +  Op\left(   R_1(g, \lambda)-R_2(g, \lambda)\right ) \\&\qquad +Op\left( \Delta {H}_{\frac{1}{2}}(g, \lambda)~X_ga(g, \lambda )-\Delta a(g, \lambda ) X_g {H}_{\frac{1}{2}}(g, \lambda)\right),
\end{align}
where $F_{g, \lambda}^{1}(\xi,u) ,F_{g, \lambda}^{2}(\xi,u) \in S^{0}(\mathbb{R}^n)$ (the term associated with  $j=1$ is $F_{g, \lambda}^{1}$ and $F_{g, \lambda}^{2}$  is the remaining terms with $j>1$  in the asymptotic expansion (\ref{realerror})). Therefore $ Op^W(F_{g, \lambda}^{1}(\xi,u)  +F_{g, \lambda}^{2}(\xi,u)) \in S_\mathcal{H}^0(\mathbb{H}^n).$  Since each symbol in the last equaity of the    expression (\ref{square}) belongs to the $S_\mathcal{H}^0(\mathbb{H}^n)$ class,  by Theorem \ref{CH80003},   the operator $[\mathcal{H}_{\frac{1}{2}}, A]$ is bounded on $L^2(\mathbb{H}^n).$

	(c)  Since $A$ is bounded self-adjoint, the spectrum of $A$, $ \sigma(A) $ is a compact subset of $\mathbb{R}.$ Since any continuous function can be approximated in the supremum norm by smooth functions, it is enough to assume that $f \in {C^2(\sigma (A))}.$ By Theorem 1.6  of Laptev-Safarov \cite{LapSaff},  by setting $A=\sqrt{\mathcal{H}}, B= A, \chi = 0, \psi = f, P_{\lambda}=\mathcal{ P}_{r^2} $, we get
	$$ | {tr}(\mathcal{P}_{r^2} f(A)\mathcal{P}_{r^2} - \mathcal{P}_{r^2} f(\mathcal{P}_{r^2}A\mathcal{P}_{r^2}) \mathcal{P}_{r^2} )|$$
	$$\leq  \frac{1}{2}\|f''\|_{\infty}N_{r_1}({r^2}) \Big({\| \mathcal{P}_{{r^2-r_1,r^2} } A\|^2+\frac{\pi ^2}{6r_1^2}\| \mathcal{P}_{r^2 - r_1}[A, \sqrt{\mathcal{H}} ]\|^2}\Big) .$$
	Dividing both sides by $tr (\mathcal{P}_{r^2})$ and setting $r_1=r^{2\alpha }, \alpha \in (0,1)$
$$\frac{| {tr}(\mathcal{P}_{r^2}f(A)\mathcal{P}_{r^2}- \mathcal{P}_{r^2} f(\mathcal{P}_{r^2} A\mathcal{P}_{r^2}) \mathcal{P}_{r^2}) |}{tr (\mathcal{P}_{r^2})} \leq C \frac{N_{r^{2\alpha }}(r)}{tr (\mathcal{P}_{r^2})}\approx r^{2\alpha-2}.$$
So	$$\frac{| {tr}(\mathcal{P}_rf(A)\mathcal{P}_r- \mathcal{P}_r f(\mathcal{P}_rA\mathcal{P}_r) \mathcal{P}_r) |}{tr (\mathcal{P}_r)} \leq C \frac{N_{r^{\alpha }}(r)}{tr (\mathcal{P}_r)} \approx r^{\alpha-1}\to 0 ~\mbox{as}~ r \to \infty $$ by part (3) and part (4) of Corollary \ref{abs}, where $N_{r_1}(r)=\displaystyle \sup_{\mu \leq r} \left( tr(\mathcal{P}_{\mu+r_1 }- \mathcal{P}_{\mu }) \right)$.
\end{proof}

The proof of Theorem \ref{sch} follows from Corollary \ref{final} and part (c) of Lemma \ref{direct}.
\section{Szeg\"o limit theorem for $\mathcal{H}_1$}\label{sec6}
Consider the  operators $\mathcal{H}_1$ and $\mathcal{H}$ as defined in Theorem \ref{bdd}. Since the operators  $e^{-t\mathcal{H}_1}$ and $e^{-t\mathcal{H}}$  are compact for $t>0$, we choose a suitable  $m \in \mathbb{N}$ such that $(\mathcal{H}_1+rI)^{-m}$ and $(\mathcal{H}+rI)^{-m}$ are trace class operators on $L^2(\mathbb{H}^n)$ for $r>0$. We observe the following facts before proving Theorem \ref{bdd}.
\begin{lemma} \label{11}
Consider the self-adjoint operators $\mathcal{H}$ and $\mathcal{H}_1$ as defined in Theorem \ref{bdd}. Then \\
(a) \begin{align*}
\bigg|\frac{tr~((\mathcal{H}_1+rI)^{-m})}{tr~((\mathcal{H}+rI)^{-m})} - 1\bigg|
\rightarrow 0  ~ \mathrm{as} ~ r \rightarrow \infty.
\end{align*}
(b) If $B$ is any bounded operator on $L^2(\mathbb{H}^n)$, then
 \begin{align*}
\bigg|\frac{tr~(B(\mathcal{H}_1+rI)^{-m})}{tr~(B(\mathcal{H}+rI)^{-m})} - 1\bigg|
\rightarrow 0  ~ \mathrm{as} ~ r \rightarrow \infty.
\end{align*}
\end{lemma}
\begin{proof}
Without loss of generality we prove the result for the positive operator $B$ by  adding a suitable constant $c>0$ which  makes the operator $B+cI$ positive.

 (a) Since $ B $ and $  (\mathcal{H}+rI)^{-1}$ are bounded and positive operators, we have
 \begin{align}
(\mathcal{H}_1+rI)=(\mathcal{H}+rI)^{\frac{1}{2}} ((\mathcal{H}+rI)^{-\frac{1}{2}} (B) (\mathcal{H}+rI)^{-\frac{1}{2}} +1     )(\mathcal{H}+rI)^{\frac{1}{2}}.
 \end{align}
 Therefore
\begin{align}\label{point}
(\mathcal{H}_1+rI)^{-m}=(\mathcal{H}+rI)^{-m} +(\mathcal{H}+rI)^{-\frac{m}{2}} ((1+K_{r})^{-m} -1)(\mathcal{H}+rI)^{-\frac{m}{2}},
\end{align}
 where $K_{r} = (\mathcal{H}+rI)^{-\frac{1}{2}} B (\mathcal{H}+rI)^{-\frac{1}{2}}$. Here $K_{r}$ is a positive operator and $\|(I+ K_{r})^{-1}\| \leq 1,$ for any $ r >0.$ Thus
  \begin{align*}
 ~\big|tr~\big ((\mathcal{H}_1+rI)^{-m})-{tr~((\mathcal{H}+rI)^{-m}\big )}\big|
= &~\big|tr~\big ((\mathcal{H}+rI)^{-\frac{m}{2}} ((1+K_{r})^{-m} -1)(\mathcal{H}+rI)^{-\frac{m}{2}}\big )\big|\\
\leq &~ tr~\big ((\mathcal{H}+rI)^{-m}\big ) \big \| ((1+K_{r})^{-m} -1)\big \|\\
\leq & ~m \big \| K_{r}\big \| tr~ \big((\mathcal{H}+rI)^{-m}\big )\\
\leq & ~m \big \| B\big \|\big \| (\mathcal{H}+rI)^{-1}\big \| tr~ \big((\mathcal{H}+rI)^{-m}\big ).
 \end{align*}
 Therefore,
 \begin{align*}
 \bigg|\frac{tr~((\mathcal{H}_1+rI)^{-m})}{tr~((\mathcal{H}+rI)^{-m})} - 1\bigg|
 \leq m\big \| B\big \|\big \| (\mathcal{H}+rI)^{-1}\big \| \rightarrow 0  ~ \mathrm{as} ~ r \rightarrow \infty.
 \end{align*}

 (b) Using (\ref{point}) we have
  \begin{align*}
 & \big|tr~\big (B(\mathcal{H}_1+rI)^{-m})-{tr~(B(\mathcal{H}+rI)^{-m}\big )}\big|\\
 &= \big|tr~\big (B(\mathcal{H}+rI)^{-\frac{m}{2}} ((1+K_{r})^{-m} -1)(\mathcal{H}+rI)^{-\frac{m}{2}}\big )\big|\\
  &= \big|tr~\big ((\mathcal{H}+rI)^{-\frac{m}{2}}B(\mathcal{H}+rI)^{-\frac{m}{2}} ((1+K_{r})^{-m} -1)\big )\big|\\
    &= \big|tr~\big (W_{r}((1+K_{r})^{-m} -1)\big )\big|\\
    &= \big|tr~\big (W_{r}^{\frac{1}{2}}((1+K_{r})^{-m} -1)W_{r}^{\frac{1}{2}}\big )\big|\\
&\leq  m \big \| B\big \|\big \| (\mathcal{H}+rI)^{-1}\big \| tr\big(B(\mathcal{H}+rI)^{-m}\big ),
 \end{align*}
 where $W_{r} = (\mathcal{H}+rI)^{-\frac{m}{2}}B(\mathcal{H}+rI)^{-\frac{m}{2}}$ is a positive and trace class operator on $L^2(\mathbb{H}^n)$.

Therefore,
\begin{align*}
\bigg|\frac{tr~(B(\mathcal{H}_1+rI)^{-m})}{tr~(B(\mathcal{H}+rI)^{-m})} - 1\bigg|\leq m\|B\|\big \| (\mathcal{H}+rI)^{-1}\big \|
\rightarrow 0  ~ \mathrm{as} ~ r \rightarrow \infty.
\end{align*}
\end{proof}

\begin{lemma} \label{22}
Let $\mathcal{H}$ and $\mathcal{H}_1$ defined as in Theorem \ref{bdd}, then
\begin{align*}
\lim_{r\to\infty}\frac{tr~(B(\mathcal{H}_1+rI)^{-m})}{tr~((\mathcal{H}_1+rI)^{-m})} = \lim_{r\to\infty} \frac{tr~(B(\mathcal{H}+rI)^{-m})}{tr~((\mathcal{H}+rI)^{-m})}.
\end{align*}
The above equality valid in the sense that if one of limits exist then the other also does and the limits are the  same.
\end{lemma}
\begin{proof}
For each $r >0$, we have
\begin{align} \label{frac}
\frac{\bigg(\frac{tr~(B(\mathcal{H}_1+rI)^{-m})}{tr~(B(\mathcal{H}+rI)^{-m})}\bigg)}{\bigg(\frac{tr~((\mathcal{H}_1+rI)^{-m})}{tr~((\mathcal{H}+rI)^{-m})}\bigg)} =  \frac{\bigg(\frac{tr~(B(\mathcal{H}_1+rI)^{-m})}{tr~((\mathcal{H}_1+rI)^{-m})}\bigg)}{\bigg(\frac{tr~((B(\mathcal{H}+rI)^{-m})}{tr~((\mathcal{H}+rI)^{-m})}\bigg)}.
\end{align}
 Since the left hand side has limit 1 (by part (b) of Lemma \ref{11}), the right hand side limit in (\ref{frac}) exists and equal to 1. Therefore if the numerator or the denominator in the fraction in the right hand side has a limit in (\ref{frac}), then the other also has a limit and they both agree. Therefore,
 $ \displaystyle\lim_{r\to\infty}\frac{tr~(B(\mathcal{H}_1+rI)^{-m})}{tr~((\mathcal{H}_1+rI)^{-m})} = \lim_{r\to\infty} \frac{tr~(B(\mathcal{H}+rI)^{-m})}{tr~((\mathcal{H}+rI)^{-m})}.$
 \end{proof}

\noindent{\bf Proof of theorem \ref{bdd}}:
Without loss of generality add a suitable constant to make the function $f$ positive. Then $f(A)$ is a positive operator.
Setting $\phi_\mathcal{H}(r)=tr(\mathcal{P}_r),~ \phi_{\mathcal{H}_1}(r)=tr(\mathcal{P}_r^{'}),$~$ \phi_{\mathcal{H},f}(r)=tr(\mathcal{P}_rf(A)\mathcal{P}_r)$ and $\phi_{{\mathcal{H}_1,f}}(r)=tr(\mathcal{P}_r^{'}f(A)\mathcal{P}_r^{'})$ we have
\begin{eqnarray*}
\lim_{r\to\infty} \frac{{tr}{(\mathcal{P}_r'f(A)\mathcal{P}_r')}}{tr~(\mathcal{P}_r')}  &=& \lim_{r\to\infty}\dfrac{\int_0^\infty \frac{\phi_{{\mathcal{H}_1,f}}(u)}{(1+\frac{u}{r})^{m+1}}du }{\int_0^\infty \frac{\phi_{{\mathcal{H}_1}}(u)}{(1+\frac{u}{r})^{m+1}} du }\\&=&\lim_{r\to\infty}\dfrac{\int_0^\infty \frac{\phi_{{\mathcal{H},f}}(u)}{(1+\frac{u}{r})^{m+1}}du }{\int_0^\infty \frac{\phi_{{\mathcal{H}}}(u)}{(1+\frac{u}{r})^{m+1}} du }\\&=&\lim_{r\to\infty} \frac{{tr}{(\mathcal{P}_rf(A)\mathcal{P}_r)}}{tr~(\mathcal{P}_r)}\\
& =&\lim_{r\to\infty}  \frac{\int_{G^{r}}f(a_{g, {\lambda}}(\xi, x))  \,d\xi\,dx \,dg\,d\mu(\lambda) }{\int_{G^{r}}  \,d\xi\,dx \,dg\,d\mu(\lambda)},
	\end{eqnarray*}
	(Assuming one limit exists)\\
	where $G^{r}=\{(g, \lambda, \xi, x)\in \mathbb{H}^n \times \mathbb{R}^*\times  \mathbb{R}^n\times \mathbb{R}^n : |\lambda |(1+|\xi| ^2+|x|^2)+V(g)\leq r  \}$  and $a(g, {\lambda})=Op^W(a_{g, {\lambda}}).$ We use Lemma \ref{22} for the middle equality and Theorem \ref{gp} (see Appendix) for the extreme left equalities. The extreme right equality follows from Lemma \ref{direct}.

 \begin{corollary}\label{cpt}
 	The Theorems \ref{mo}, \ref{sch} and \ref{bdd} also hold under the compact perturbation of the pseudo-differential operator $A$.
 \end{corollary}
\begin{proof}
To prove the above result, enough to show $\displaystyle
\lim_{r \rightarrow \infty} \frac{tr (\mathcal{P}_r A^n \mathcal{P}_r )}{tr~ (\mathcal{P}_r)} = \displaystyle
\lim_{r \rightarrow \infty} \frac{tr (\mathcal{P}_r (A+K)^n \mathcal{P}_r )}{tr~ (\mathcal{P}_r)} $ for any compact operator $K$ on $L^2(\mathbb{H}^n)$. Notice that $(A+K)^n = A^n$+ terms with factor $A^{p}K^{n-p} $ or $K^{p}A^{n-p}$  where $ p \in \{1,2, \cdots ,n\}$. Since the class of compact operators form a two sided ideal of the class of bounded operators, $(A + K)^n = A^n+$ a compact operator. We are done if we can prove that for a compact  operator $T$,  $\displaystyle \lim_{r \rightarrow \infty} \frac{tr (\mathcal{P}_r T \mathcal{P}_r )}{tr~ (\mathcal{P}_r)}=0.$
Since T is a compact operator, for given $ \epsilon> 0$ there exist a finite rank
 operator $T_k$ such that $\|T_k - T\| \rightarrow 0 ~ \text{as } k \rightarrow \infty$. Then
 $
\left |\dfrac{tr (\mathcal{P}_r T \mathcal{P}_r )- tr (\mathcal{P}_r T_{k} \mathcal{P}_r )}{tr (\mathcal{P}_r)} \right | \leq \|T-T_k\| ~  \rightarrow 0 ~ \text{as } k \rightarrow \infty.
 $
 Therefore for given $ \epsilon> 0$ there exist $N_0\in \mathbb{N}$ such that $\left |\dfrac{tr (\mathcal{P}_r T \mathcal{P}_r )- tr (\mathcal{P}_r T_{k} \mathcal{P}_r )}{tr (\mathcal{P}_r)} \right |< \dfrac{\epsilon}{2}$ for $k \geq  N_0$. Further,  $\bigg |\frac{tr (\mathcal{P}_r T_{N_0} \mathcal{P}_r )}{tr~ (\mathcal{P}_r)}\bigg |  \rightarrow 0 ~ \text{as } r \rightarrow \infty$ i.e, for given $\epsilon >0, ~~\exists ~N_1 \in \mathbb{N}$ such that  $\bigg |\frac{tr (\mathcal{P}_r T_{N_0} \mathcal{P}_r )}{tr~ (\mathcal{P}_r)}\bigg | < \dfrac{\epsilon}{2},~~ \forall  ~ r> N_1.$ Thus
$
\bigg |\frac{tr (\mathcal{P}_r T \mathcal{P}_r )}{tr~ (\mathcal{P}_r)} \bigg |
\leq  \bigg |\frac{tr (\mathcal{P}_r T_{N_0} \mathcal{P}_r )}{tr~ (\mathcal{P}_r)}\bigg |+ \|T-T_{N_0}\|
< \epsilon \quad \forall ~r \geq  N_1.
 $
\end{proof}
\begin{remark}\label{final remark}
     The proof of part (c) of Lemma \ref{direct} can also be achieved for $\kappa\in (0,1)$ proving the boundedness of the operators $[A, V]$ and $[A, \mathcal{L}]$ on $L^2(\mathbb{H}^n)$. Now for any  $h  \in L^2(\mathbb{H}^n)$, we have 
	$$	[A, V] h(g)=(AV-VA)h(g)=\int_{\mathbb{H}^n} K_3(g, g_1)h(g_1) \, d g_1,$$
		where  
		\begin{align}\label{CH80004}
			K_3(g, g_1)=	\big (V(g_1)- V(g)\big )\int_{\mathbb{R}^*} {tr}\big(\pi_{\lambda}(g)~a(g, \lambda)\pi_{\lambda}^*(g_1)\big ) \, d\mu(\lambda).
		\end{align}We note that 
\begin{align}\label{CH8note}\nonumber
a(g, \lambda)&=\widehat{k_g}(\lambda)\\\nonumber
&=\int_{\mathbb{H}^n} k_g(g_1)\pi_{\lambda}^*(g_1)dg_1\\\nonumber
&=\int_{\mathbb{H}^n} k_g(g_1)	(I-T^2)^{N}(1+\lambda^2)^{-N}\pi_{\lambda}^*(g_1)dg_1\\\nonumber
&=(-1)^N(1+\lambda^2)^{-N}\int_{\mathbb{H}^n} (I-T^2)^{N} k_g(g_1)	\pi_{\lambda}^*(g_1)dg_1\\				\nonumber
&=(-1)^N(1+\lambda^2)^{-N}\widehat{\left((I-T^2)^{4N} k_g\right)}(\lambda)\\
&=(-1)^N(1+\lambda^2)^{-N} \sum_{\beta\leq 2N} \pi(T)^\beta a(g, \lambda).
\end{align}
Then, using the identity   (\ref{CH8note}), we have
\begin{align}\label{convolution}\nonumber
&|{tr}\big(\pi_{\lambda}(g)~a(g, \lambda)\pi_{\lambda}^*(g_1)\big )|=|{tr}\left(\pi_{\lambda}(g_1^{-1}g)~a(g, \lambda)\right )|\\\nonumber
&=\bigg|{tr}\bigg(\pi_{\lambda}(g_1^{-1}g)\left(\pi_{\lambda}(I-\mathcal{L})+V(g_1^{-1}g)\right)^{-4N}\left(\pi_{\lambda}(I-\mathcal{L})+V(g_1^{-1}g)\right)^{4N}\bigg.\\\nonumber
&\bigg.\qquad \times (1+\lambda^2)^{-4N}  \sum_{\beta\leq 4N} \pi(T)^\beta a(g, \lambda)\bigg )\bigg|\\\nonumber
&\leq \left |{tr}\big(\pi_{\lambda}(I-\mathcal{L})+V(g_1^{-1}g)\big)^{-4N}\right |  \left \|\left(\pi_{\lambda}(I-\mathcal{L})+V(g_1^{-1}g)\right)^{4N}(1+\lambda^2)^{-4N}\right\|_{op} \\\nonumber&\qquad \times \bigg \| \sum_{\beta\leq 4N} \pi(T)^\beta a(g, \lambda)\bigg \|_{op}\\
&\leq C \left |{tr}\big(\pi_{\lambda}(I-\mathcal{L})+V(g_1^{-1}g)\big)^{-4N}\right|,
\end{align}
where the second term is bounded by Theorem \ref{CH80003} and the last term is bounded by $\|a\|$.
Now  we have	\begin{eqnarray*}  
tr ~(1+|\lambda|H+V(g_1^{-1}g))^{-4N}&=&\sum_\alpha\langle(1+|\lambda|H+V(g_1^{-1}g))^{-4N}\Phi_\alpha,\Phi_\alpha\rangle \\\nonumber &=& \sum_{\alpha} \frac{1}{(1+|\lambda| (2|\alpha|+n)+V(g_1^{-1}g))^{4N}}\\\nonumber
\end{eqnarray*}
and consequently
\begin{align}\label{CH10001}\nonumber
&\int_{\mathbb{R}^*} {tr}|\big(\pi_{\lambda}(g)~a(g, \lambda)\pi_{\lambda}^*(g_1)\big ) |\, d\mu(\lambda) \\\nonumber
&\leq C \sum_{\alpha}  	\int_{0}^\infty  \frac{\lambda^n}{(1+\lambda (2|\alpha|+n)+V(g_1^{-1}g))^{4N}} \, d\lambda\\\nonumber
&= C \sum_{\alpha}  \frac{1}{(2|\alpha|+n)^{n+1}  }	\int_{1+V(g_1g^{-1})}^\infty  \frac{(u-1-V(g_1^{-1}g))^n}{u^{4N}} \, du\\ 
&= C_2 \frac{1}{(1+V(g_1^{-1}g))^{4N-n-1}}.
\end{align}
Thus  from (\ref{CH80004}), (\ref{convolution}), and (\ref{CH10001}),  we get $$ |K_3(g, g_1)|\leq \frac{c|(V(g_1)- V(g)|}{(1+V(g_1^{-1}g))^{4N-n-1}}.$$ For large $|g|, |g_1|$, using  the trangle inequality for the homogeneous norm and the fact that $\kappa\in (0,1)$, we have
\begin{align*}
	\|[A, V] h\|_2^{2}\leq &\int_{\mathbb{H}^n} \bigg |  \int_{\mathbb{H}^n} K_3(g, g_1)h(g_1) \, d g_1 \bigg |^2 \, dg
	\\&\leq c_N \int_{\mathbb{H}^n} \bigg (  \int_{\mathbb{H}^n}  \bigg |  \frac{\big |  |g_1|^k-|g|^k\big | h(g_1) }{(1+\big |  g_1^{-1}g\big |^k)^{4N-n-1}} \bigg |  \, d g_1 \bigg )^2 \, d g\\
		\\&\leq c_N \int_{\mathbb{H}^n} \bigg (  \int_{\mathbb{H}^n}  \bigg |  \frac{\big |  |g_1|-|g|\big | h(g_1) }{(1+\big |  g_1^{-1}g\big |^k)^{4N-n-1}} \bigg |  \, d g_1 \bigg )^2 \, d g\\
	\\&\leq c_N \int_{\mathbb{H}^n} \bigg (  \int_{\mathbb{H}^n}  \bigg |  \frac{\big |  g_1^{-1}g\big | h(g_1) }{(1+\big |  g_1^{-1}g\big |^k)^{4N-n-1}} \bigg |  \, d g_1 \bigg )^2 \, d g\\
&	\leq c_N \int_{\mathbb{H}^n} \bigg (  \int_{\mathbb{H}^n}  \bigg |  \frac{   h(g_1) }{(1+\big | g_1^{-1}g\big |^k)^{4N-n-2}} \bigg |  \, d g_1 \bigg )^2 \, d g
	= \| |h| \ast K\|_2^2,
	\end{align*}where $K(g)= \dfrac{1}{(1+|g|^k)^{ 4N-n-2}}$. Since  for a sufficiently large $N\in\mathbb{N}$, $K \in L^1({\mathbb{H}^n) }$,  an application of Minkowski's inequality gives $\|[A, V]h\|_2 \leq C \|K\|_1 \|h\|_2.$

If $|g|$ and $|g_1|$ are lying in some compact set $\mathcal{K}\subset \mathbb{R}$ then $ \int_{\mathcal{K}} \left |  \int_{\mathcal{K}} K_3(g, g_1)h(g_1) \, d g_1 \right |^2 \, dg\leq C_{\mathcal{K}}\|h\|_2$. If $|g| $ (or $|g_1|$) lies in $K$ and $|g_1|$ (or $|g|$) is large then an application of Cauchy-Schwarz inequality gives $\|[A, V]h\|_2 \leq \|h\|_2\int_{\mathcal{K}}  \int_{|g|} |K_3(g, g_1)|^2 \, d g_1   \, dg\leq C_{\mathcal{K}}\|h\|_2$.

For $\kappa\in(0, 1)$ the operator $[ V, A]$ is bounded.  The boundedness of the operator $[\mathcal{L}, A]$ will imply boundedness of the operator $[\mathcal{H}, A]$ on $L^2(\mathbb{H}^n)$ as $[T^2, A]=0$. Using the identity (\ref{CH8note}),  we get
\begin{align*}
A\mathcal{L} h(g)&=\int_{\mathbb{R}^*} {tr}\left (\pi_{\lambda}^*(g)a(g, \lambda)\widehat{\mathcal{L} h}(\lambda)\right ) \,  d\mu(\lambda)\\
&=\int_{\mathbb{R}^*} {tr}\left ( \pi_{\lambda}^*(g) (1+\lambda^2)^{-2N}   \sum_{\beta\leq 4N} \pi(T)^\beta a(g, \lambda) ~ \hat{h}(\lambda) |\lambda|{H}\right ) \,  d\mu(\lambda)\\
&=\int_{\mathbb{R}^*} {tr}\bigg (  \left(\pi_{\lambda}^*(g) \pi_{\lambda}(I-\mathcal{L})+V(g)\right)^{-4N+1}\left(\pi_{\lambda}(I-\mathcal{L})+V(g)\right)^{4N}(1+\lambda^2)^{-2N}\\&\quad \quad \times  \sum_{\beta\leq 4N} \pi(T)^\beta a(g, \lambda) ~ \hat{h}(\lambda) |\lambda|{H}  \left( \pi_{\lambda}(I-\mathcal{L})+V(g)\right)^{-1}\bigg ) \,  d\mu(\lambda).
\end{align*}
Arguing as in (\ref{convolution}) and (\ref{CH10001}), we obtain 
\begin{align*}
|A\mathcal{L} h(g)|&\leq \int_{\mathbb{R}^*} {tr}\left (  \left( \pi_{\lambda}(I-\mathcal{L})+V(g)\right)^{-4N+1}\right) \|\hat{h}(\lambda) \|_{op} \,  d\mu(\lambda)\\
&\leq \int_{\mathbb{R}^*} {tr}\left (  \left( \pi_{\lambda}(I-\mathcal{L})+V(g)\right)^{-4N+1}\right) \|\hat{h}(\lambda) \|_{S_2} \,  d\mu(\lambda)\\
&\leq \left( \int_{\mathbb{R}^*} \left| {tr}\left (  \left( \pi_{\lambda}(I-\mathcal{L})+V(g)\right)^{-4N+1}\right) \right|^2 \,  d\mu(\lambda)\right)^{\frac{1}{2}}  \left( \int_{\mathbb{R}^*}  \|\hat{h}(\lambda)  \|_{S_2}^2\, d\mu(\lambda)\right)^{\frac{1}{2}}\\
& \leq C \|h\|_2 (1+V(g))^{\frac{-8N+n+3}{2}} .
\end{align*}
For sufficiently large $N$, an application of Cauchy-Schwarz inequality gives
$\| A\mathcal{L}h \|_2 \leq M_1 \|h\|_2.$
Further,
\begin{align*}
\mathcal{L}A h(g)&=\int_{\mathbb{R}^*} {tr}\left (\pi_{\lambda}^*(g)\widehat{\mathcal{L}A h}(\lambda)\right ) \,  d\mu(\lambda)\\\nonumber
&=\int_{\mathbb{R}^*} {tr}\left (\pi_{\lambda}^*(g)  \widehat{Ah}(\lambda) |\lambda|{H}\right ) \,  d\mu(\lambda)\\\nonumber
&=\int_{\mathbb{R}^*} \int_{\mathbb{H}^n} Ah(g_1) ~{tr}\bigg ( \pi_{\lambda}^*(g)  \pi_{\lambda}(g_1) |\lambda| {H}\bigg ) \,  d\mu(\lambda) \, dg_1\\\nonumber
&=\int_{\mathbb{R}^*} \int_{\mathbb{H}^n}\left( \int_{\mathbb{R}^*} {tr}\left (\pi_{\lambda_1}^*(g_1)a(g_1, \lambda_1)\widehat{ h}(\lambda_1 )\right ) \,  d\mu(\lambda_1) \right )  ~{tr}\bigg ( \pi_{\lambda}^*(g)  \pi_{\lambda}(g_1) |\lambda| {H}\bigg ) \,  d\mu(\lambda) \, dg_1.
\end{align*}
 Proceeding as   in (\ref{CH10001}), we get
$$
\displaystyle\int_{\mathbb{R}^*} \left |{tr}\left (\pi_{\lambda_1}^*(g_1)a(g_1, \lambda_1)\widehat{ h}(\lambda_1 )\right ) \right| \,  d\mu(\lambda_1)\leq C\|h\|_2 (1+V(g_1))^{\frac{-8N+n+1}{2}} .
$$
Moreover,   a similar way as in   (\ref{CH10001}), we obtain  \begin{eqnarray*}&&\int_{\mathbb{R}^*} \int_{\mathbb{H}^n}\left| {tr}\left ( \pi_{\lambda}^*(g)  \pi_{\lambda}(g_1) |\lambda| {H} \right )\right|  (1+V(g_1))^{\frac{-8N+n+1}{2}} \,  d\mu(\lambda) \, dg_1\\
	&\leq & \int_{\mathbb{H}^n}  (1+V(g_1))^{\frac{-8N+n+1}{2}}  dg_1 \int_{\mathbb{R}^*} \left| {tr}\left (   |\lambda| {H}~ \pi_{\lambda}^*(g) \right  )\right|   \,  d\mu(\lambda) \\
&\leq &  \int_{\mathbb{R}^*}tr \left ( (1+\lambda^2)^{-4N}  \sum_{\beta\leq 8N} \pi(T)^\beta |\lambda| {H}\right ) \,  d\mu(\lambda) \\
&\leq& C (1+V(g))^{-8N+n}  
\end{eqnarray*}
Therefore $\|[A, \mathcal{L}]h\|_2 \leq M  \|h\|_2$ and so the operator $[A, \mathcal{H}]$ is bounded on $L^2(\mathbb{H}^n)$.

Now setting $A=\mathcal{H}, B= A, \chi = 0, \psi = f, P_\lambda = \pi_{r } $ in Theorem 1.6  of Laptev-Safarov \cite{LapSaff}, we get
		$$ | {tr}(\mathcal{P}_rf(A)\mathcal{P}_r - \mathcal{P}_r f(\mathcal{P}_rA\mathcal{P}_r) \mathcal{P}_r )|$$
	$$\leq  \frac{1}{2}\|f''\|_{\infty}N_{r}(r) \Big({\| \pi_{r } A\|^2+\frac{\pi ^2}{6r^2}\|\pi_{r - r}[A, \mathcal{H} ]\|^2}\Big) .$$
	Dividing both sides by $tr (\mathcal{P}_r)$ and setting $r=r^{\alpha }, \alpha \in (0,1)$ and using the boundness of $A, [A, \mathcal{H} ]$ we have
	$$\frac{| {tr}(\mathcal{P}_rf(A)\mathcal{P}_r- \mathcal{P}_r f(\mathcal{P}_rA\mathcal{P}_r) \mathcal{P}_r) |}{tr (\mathcal{P}_r)} \leq C \frac{N_{r^{\alpha }}(r)}{tr (\mathcal{P}_r)} \to 0 ~\mbox{as}~ r \to \infty $$ by (4) of Corollary \ref{abs}, where $N_{r}(r)=\displaystyle \sup_{\mu \leq r} \left( tr(\pi_{\mu }- \pi_{\mu-r }) \right)$.

\end{remark}
\section*{Acknowledgments}
The first author wishes to thank the Ministry of Human Resource Development, India for the  research fellowship and Indian Institute of Technology Guwahati, India for the support provided during the period of this work.

\section{Appendix}

We collect few definitions and theorems of Grishin-Poedintseva \cite{gri}, that we use in our paper for the reader's convenience.
\begin{definition}
 Let $\phi$ be a positive function on the half line $[0,\infty)$. Let
$$
S = \{ \alpha : \exists M,R ~ \mathrm{with} ~ \phi(tr) \leq M t^\alpha\phi(r), ~ \mathrm{for ~ all ~ }  t \geq 1, r \geq R \}
$$
and
$$
G = \{ \beta : \exists M,R ~ \mathrm{with} ~ \phi(tr) \geq M t^\beta\phi(r), ~ \mathrm{for ~ all ~ }  t \geq 1, r \geq R \}
$$
Then the numbers {\bf $\alpha(\phi):=\inf S$} and
{\bf $\beta(\phi):=\sup G$}  are called the {\bf upper} and {\bf lower Matushevskaya index}
 of $\phi$ respectively.

\end{definition}
\begin{theorem}(\cite{gri},Theorem 2)\label{r}

Let $m>-1$. Assume that $\varphi$ is positive measurable function on $[0,\infty)$ that
does not vanish identically in any neighborhood of infinity. Let
$\Phi(r)=\displaystyle\int_0^\infty \dfrac{\varphi(rt)}{(1+t)^{m+1}}dt$
be finite.
Then the functions $\varphi$ and $\Phi$ have same growth at infinity if and only if
$\beta(\varphi)>-1$ and $\alpha(\varphi)<m$.
\end{theorem}
\begin{definition}
A function $\varphi$ is said to be multiplicatively continuous at infinity if it satisfies
$\displaystyle\lim_{\substack{r\rightarrow \infty\\\tau\rightarrow 1}}\dfrac{\varphi(\tau r)}{\varphi(r)}=1.$
\end{definition}

\begin{theorem}(\cite{gri},Theorem 8)\label{gp}
 Let $\varphi$ and $\psi$ be positive functions on $[0,\infty)$ satisfying the following conditions:
\begin{enumerate}
 \item the functions $\varphi$ and $\psi$ do not vanish identically in any neighborhood of infinity;
\item the function $\varphi$ is multiplicatively continuous at infinity and $\beta(\varphi)>-1$;
\item the function $\psi$ is increasing;
\item at least one of the inequalities $\alpha(\varphi)<m$ and $\alpha(\psi)<m$ holds, where $m>-1$;
\item the functions
$$
\Phi(r)=\int_0^\infty \dfrac{\varphi(ru)}{(1+u)^{m+1}}du ~~and~~\Psi(r)=\int_0^\infty \dfrac{\psi(ru)}{(1+u)^{m+1}}du
$$
are finite and $\displaystyle\lim_{r\to \infty}\dfrac{\Psi(r)}{\Phi(r)}=1$ then $\displaystyle\lim_{r\to \infty}\dfrac{\psi(r)}{\varphi(r)}=1$.
\end{enumerate}
 \end{theorem}


\begin{thebibliography}{99}
	\normalsize
\baselineskip=17pt


\bibitem{atiyah}  M. Atiyah, R.  Bott, and V. K. Patodi,  On the heat equation and the index theorem, \emph{Invent. Math.} 19, 279-330 (1973). 

\bibitem{BKG} H. Bahouri, C. Fermanian-Kammerer, and  I. Gallagher,  Phase-space analysis and pseudodifferential calculus on the Heisenberg group, \emph{Ast{\'e}risque} (342),    (2012).


\bibitem{ruz2014} V. Fischer  and  M. Ruzhansky, A pseudo-differential calculus on the Heisenberg group, \emph{C. R. Math. Acad. Sci. Paris} 352(3),  197-204 (2014).



\bibitem{ruz14} V. Fischer and  M. Ruzhansky, \emph{Quantization on nilpotent Lie groups},  Progress in Mathematics, vol. 314, Birkh\"auser, Basel (2016).


\bibitem{gri}  A. F. Grishin and I. V. Poedintseva, Towards the Keldysh Tauberian theorem, \emph{J. Math. Sci.} 134(4), 2272-2287 (2006).






\bibitem{gui}   V. Guillemin, \emph{Some classical theorems in spectral theory revisited}, In: Seminar on singularity of solutions of Linear Partial  differential equations, pp.  219-259, Princeton Univ. Press,  Princeton  (1979).

\bibitem{hor}   L. H\"{o}rmander,  \emph{The Analysis of Linear Partial Differential Operators IV},  Springer-Verlag, Berlin (1985).








\bibitem{JZ} A. J. E. M. Janssen and  S. Zelditch,  Szeg\"o limit theorems for the harmonic oscillator, \emph{Trans. Amer. Math. Soc.}   280(2), 563-587 (1983).


\bibitem{kel51} M. V. Keldysh,  On a Tauberian theorem, \emph{Trudy Mat. Inst. Steklov.} 38,  77-86 (1951).



\bibitem{kr} B. Kr\"otz, Bernhard, S. Thangavelu, and Y.  Xu,  The heat kernel transform for the Heisenberg group, \emph{J. Funct. Anal.} 225(2),  301-336 (2005).



\bibitem{tan} H. Kumano-go and  K.  Taniguchi,
Oscillatory integrals of symbols of pseudo-differential operators on $\mathbb{R}^n$ and operators of Fredholm type,
\emph{Proc. Japan Acad.} 49, 397-402 (1973).



\bibitem{lap91} A. Laptev and  Yu. Safarov,  Error estimate in the generalized Szeg\"o theorem, \emph{Équations aux Dérivées Partielles}, XV-1-XV-7, Saint-Jean-De-Monts (1991). 

\bibitem{LapSaff}  A. Laptev and  Yu. Safarov,  Szeg\"o type limit theorems, \emph{J. Funct. Anal.}  138(2),  544-559 (1996).

\bibitem{ler} N. Lerner,  \emph{Metrics on the phase space and non-selfadjoint pseudo-differential operators}, In:  Pseudo-Differential Operators,  Theory and Applications, volume 3, Birkh\"auser Verlag, Basel (2010).




\bibitem{Jiman}  L. Peng and  J. Zhao, Weyl transforms associated with the Heisenberg group, \emph{Bull. Sci. Math.} 132(1), 78-86 (2008).
 
\bibitem{rob83}  D. Robert, Remarks on the paper of S.  Zelditch: ``Szeg\"o limit theorems in quantum mechanics'', \emph{J. Funct. Anal.} 53(3),  304-308 (1983).

 
\bibitem{si} B. Simon, Schr\"odinger operators with purely discrete spectrum, \emph{Methods Funct. Anal. Topology} 15(1), 61-66 (2009).



\bibitem{sm}  J. Swain  and   M. Krishna,  Szeg\"o limit theorem on the lattice, \emph{J. Pseudo-Differ. Oper. Appl.} 10(2), 489-503 (2019).




\bibitem{tha98} S. Thangavelu, \emph{Harmonic Analysis on the Heisenberg Group},  Progress in Mathematics, vol.  159,  Birkh\"auser, Boston  (1998).




\bibitem{wid41}   D. Widder,  \emph{The Laplace transform}, Princeton Mathematical Series, vol. 6, Princeton Univ. Press,  Princeton (1941).

\bibitem{wim}  H. Widom, Eigenvalue distribution theorems for certain homogeneous spaces, \emph{J. Funct. Anal.}  32(2),  139-147 (1979).


\bibitem{zel83}  S. Zelditch, Szeg\"o limit theorems in quantum mechanics, \emph{J. Funct. Anal.}  50(1),   67-80 (1983).

 

\end{thebibliography}
\end{document}